\documentclass{amsart}
\usepackage{amssymb,amsmath}

\newcommand{\RR}{{\mathbb R}}

\newcommand{\e}{\varepsilon}

\newcommand{\be}{\beta}
\newcommand{\Lap}{\Delta}
\newcommand{\de}{\delta}
\newcommand{\del}{\partial}
\newcommand{\om}{\omega}

\newcommand{{\loc}}{{\ell\mathrm oc}}

\def\meanint{{\diagup\hskip -.42cm\int}}

\newtheorem{theorem}{Theorem}
\newtheorem{lemma}{Lemma}
\newtheorem{proposition}{Proposition}

\newtheorem{cor}{Corollary}

\begin{document}

\title{On the Fundamental Solution of an Elliptic Equation in Nondivergence Form}
\author{Vladimir Maz'ya \and Robert McOwen}
\address{Vladimir Maz'ya at Link\"oping University, Ohio State University, University of Liverpool \and Robert McOwen at Northeastern University}
\date{June 23, 2008}

\maketitle

\begin{abstract}
We consider the existence and asymptotics for the fundamental solution of an elliptic operator in nondivergence  form, ${\mathcal L}(x,\del_x)=a_{ij}(x)\del_i\del_i$, for $n\geq 3$. We assume that the coefficients have modulus of continuity  satisfying the square Dini condition.  For fixed $y$, we construct a solution of ${\mathcal L}Z_y(x)=0$ for $0<|x-y|<\e$ with explicit leading order term which is $O(|x-y|^{2-n}e^{I(x,y)})$ as $x\to y$, where $I(x,y)$  is given by an integral and plays an important role for the fundamental solution: if $I(x,y)$ approaches a finite limit as $x\to y$, then we can solve ${\mathcal L}(x,\del_x)F(x,y)=\de(x-y)$, and $F(x,y)$ is asymptotic as $x\to y$ to the fundamental solution for the constant coefficient operator ${\mathcal L}(y,\del_x)$. On the other hand, if $I(x,y)\to -\infty$ as $x\to y$ then the solution $Z_y(x)$ violates the ``extended maximum principle'' of Gilbarg \& Serrin \cite{GS} and is a distributional solution of ${\mathcal L}(x,\del_x)Z_y(x)=0$ for $|x-y|<\e$ although $Z_y$ is not even  bounded as $x\to y$.



\end{abstract}

\section{Introduction}
\smallskip\noindent
{\bf Background.}
Consider an elliptic operator in nondivergence form
\begin{equation}
{\mathcal L}(x,\del_x)\,u(x)=\,a_{ij}(x)\,\del_i\del_ju(x),
\label{eq:L}
\end{equation}
where $\del_i=\del/\del x_i$ and we have used the summation convention for repeated indices. The coefficients $a_{ij}=a_{ji}$ are real-valued functions defined on $\RR^n$ for $n\geq 3$, and we denote   the symmetric and positive definite matrix $(a_{ij}(x))$ by ${\bf A}_x$. (The case $n=2$ can be treated with a similar analysis, but additional complications arise which we have chosen to avoid here.) A {\it fundamental solution} for ${\mathcal L}$ in an open set $U$ is a function $F(x,y)$ satisfying $F(x,\cdot)\in L^1_\loc(U)$ and
\begin{equation}
-{\mathcal L}(x,\del_x)F(x,y)=\delta(x-y) \quad\hbox{for}\ x,y\in U
\label{def:FS}
\end{equation}
 in a distributional sense that needs to be made clear; for this some regularity of the coefficients will be required. If $F(x,y)$ satisfying (\ref{def:FS}) exists, then the operator $Q$ defined by
\begin{equation}
Q\phi(x)=-\int_U F(x,y)\,\phi(y)\,dy
\label{def:Q}
\end{equation}
provides a right-inverse for ${\mathcal L}$ on $C_0(U)$, the space of continuous functions with compact support in $U$.
 
 In the ``classical'' case that the coefficient functions are $\lambda$-H\"older continuous in a bounded domain $U$ for some $\lambda\in (0,1)$, it is well-known (cf.\cite{M}) that a fundamental solution exists in $U$ and is asymptotic (as $x\to y$) to the fundamental solution for the constant coefficient operator obtained by freezing the coefficients $a_{ij}$ at $y$: for $n\geq 3$ this means
\begin{equation}
 F(x,y)=\tilde F_y(x-y)(1+ H(x,y)),
\label{eq:F-asym-Holder}
\end{equation}
where, letting  $\langle , \rangle$ denote the inner product in $\RR^n$,
\begin{equation}
\tilde F_y(x)=\frac{\langle {\bf A}_y^{-1}x,x\rangle^{\frac{2-n}{2}}}{(n-2)\,|S^{n-1}|\,\sqrt{\hbox{det}{\bf A}_y}}
\label{def:F_y}
\end{equation}
 is the fundamental solution for the constant coefficient operator
${\mathcal L}(y,\del_x)\!=\!a_{ij}(y)\del_i\del_j,
$
and the remainder term $H(x,y)$ in (\ref{eq:F-asym-Holder}) satisfies
\begin{equation}
 |H(x,y)|+r|D_x H(x,y)|+r^2|D_x^2 H(x,y)|\leq c\, r^{\lambda}\quad\hbox{as}\ \ r=|x-y|\to 0,
\label{eq:H-est}
\end{equation}
for all $y$ in a compact subset of $U$. This fundamental solution may be used (cf.\cite{M}) to show the classical regularity result: if $u\in C^2(U)$ and ${\mathcal L}u$ is $\lambda$-H\"older continuous in $U$, then
$\del_i\del_j u$ is $\lambda$-H\"older continuous in $U$. 

The H\"older continuity may be generalized by assuming the coefficients have a weaker modulus of continuity, i.e.\ $a_{ij}\in C^{\om}(U)$ where $\om(r)$ is a continuous, nondecreasing function for $0\leq r< 1$ satisfying $\om(0)=0$, and
\begin{equation}
C^{\om}(U)=\{ f\in C(U): |f(x)-f(y)|\leq c\,\om(|x-y|) \quad\hbox{for}\ x,y\in U\}.
\label{def:C_om}
\end{equation}
If $\om$ satisfies  the Dini condition at zero, i.e.\ $\int_0^1\om(t)t^{-1}\,dt<\infty$, then we say that 
the coefficients are {\it Dini continuous}. In this case, there are regularity results analogous to the case of H\"older continuity  (cf.\ Proposition 1.14 in Chapter 3 of \cite{T});
however, we could not find in the literature an asymptotic description of the fundamental solution such as (\ref{eq:F-asym-Holder}) with estimates on the second-order derivatives $D^2_x H(x,y)$.  

Dini continuity is also essential for the ``extended maximum principle'' of Gilbarg \& Serrin \cite{GS} to hold: a $C^2$-solution of 
\[
{\mathcal L}u\geq 0\ \hbox{for}\  0<|x|\leq r_0
\]
 with 
 \[
 u(x)=o(|x|^{2-n})\ \hbox{ as} |x|\to 0
 \]
  must satisfy 
 \[
 u(x)<M:=\max\{u(y):|y|=r_0\}\ \hbox{for}\  0<|x|< r_0,
 \]
 and $\limsup_{|x|\to 0}u(x)<M$. In fact, they give an example (which we will discuss in Section 2) in which the coefficients are not Dini continuous and the extended maximum principle fails.

The above regularity assumptions (H\"older or Dini continuity) on the coefficients are required to study the behavior of the fundamental solution as a function of $x$ (for fixed $y$). If we instead fix $x$ and consider the behavior in $y$, then regularity of the coefficients $a_{ij}$ is not required; however, we cannot expect to achieve as precise an asymptotic description as (\ref{eq:F-asym-Holder}).  This is most conveniently described in terms of the {\it Green's function} for (\ref{eq:L}) on a smooth, bounded domain $U$, which may be defined (as in \cite{E}) to be
$G(x,\cdot)\in L^1_\loc(U)$ satisfying
\begin{equation}
\phi(x)=-\int_{U}G(x,y)\,{\mathcal L}(y,\del_y)\,\phi(y)\,dy \ \hbox{for}\ \phi\in C^2(\overline{U})\ \hbox{with}\ \phi=0\ \hbox{on}\ \del U.
\label{def:Green's}
\end{equation}
Notice that (\ref{def:Green's}) can be expressed formally as $-{\mathcal L}^*(y,\del_y)G(x,y)=\de(x-y)$  and implies that $Q\phi(x)=-\int_U G(x,y)\,\phi(y)\,dy$ defines a left-inverse for ${\mathcal L}$ on $C_0^2(U)$.
When the $a_{ij}$ are measurable, bounded, and uniformly elliptic in $U$, then the Green's function is known to exist and Fabes \& Strook \cite{FS} showed that $G(x,\cdot)\in L^{q}(U)$ for some $q>n/(n-1)$, while Bauman \cite{B1}, \cite{B2}, \cite{B3}, and Escauriaza  \cite{E} obtained pointwise estimates on $G(x,y)$ as $y\to x$ in terms of a nonnegative ``adjoint solution'' $W(y)$ which satisfies ${\mathcal L}^*(y,\del_y)W(y)=0$ in $U$. However, our paper is not concerned with such general coefficients, and for us a Green's function will also be a fundamental solution in the sense of (\ref{def:FS}).

\smallskip\noindent
{\bf Our results.}
In this paper, we allow the coefficients $a_{ij}$ to be less regular than Dini continuous, and we want to study the solutions of ${\mathcal L}u(x)=0$ with an isolated singularity at $x=y$, as well as the existence and asymptotics of  a fundamental solution $F(x,y)$ satisfying (\ref{def:FS}) in an appropriate distributional sense. We assume that the coefficients have modulus of continuity $\om$ satisfying the ``square-Dini condition''
\begin{equation}
\int_0^1 \om^2(t)\,\frac{dt }{t}< \infty.
\label{eq:om1}
\end{equation}
Condition (\ref{eq:om1}) has been encountered by other authors in different contexts: cf.\ \cite{FJK}, \cite{F}, \cite{KP}, \cite{SZ}.

To construct our solution of (\ref{def:FS}), we first fix $y$ and seek a solution of
\begin{equation}
{\mathcal L}(x,\del_x)Z_y(x)=0 \quad\hbox{for}\  x\in B_{\e}(y)\backslash\{y\},
\label{eq:LZ_y=0}
\end{equation}
where $B_{\e}(y)=\{x:|x-y|<\e\}$ for $\e$ sufficiently small, and
$Z_y(x)$ has the appropriate singularity as $x\to y$. Assuming that the modulus of continuity at $y$ satisfies (\ref{eq:om1}), we shall construct a solution of (\ref{eq:LZ_y=0}) with the asymptotic description 
\begin{equation}
Z_y(x)\sim \langle{\bf A}_y^{-1}(x-y),(x-y)\rangle^{\frac{2-n}{2}}\,
e^{I(x,y)} \quad\hbox{as}\ x\to y,
\label{eq:Zy-asym}
\end{equation}
where the factor $e^{I(x,y)}$ adjusts for lack of regularity in the coefficients: if the $a_{ij}$ are H\"older continuous, then we can take $I(x,y)\equiv 0$ and $c_y\,Z_y(x)$ is asymptotic to $\tilde F_y(x-y)$ as $x\to y$. In general, however, we find that
\begin{equation}
I(x,y)=I_y\left(\sqrt{\langle {\bf A}_y^{-1}(x-y),(x-y)\rangle}\right),
\label{eq:I(x,y)}
\end{equation}
where $I_y(r)$ is given by
\begin{equation}
\frac{1}{|S^{n-1}|}\int_{r<|z-y|<\e}
\left({\rm tr}({\bf A}_z{\bf A}_y^{-1})-
n\,\frac{\langle {\bf A}_z {\bf A}_y^{-1/2}(z-y), {\bf A}_y^{-1/2}(z-y)\rangle}{|z-y|^2}\right)\,\frac{d z}{|z-y|^{n}}
\label{def:I_y}
\end{equation}
with ${\rm tr}$ denoting matrix trace. As $r\to 0$, $I_y(r)$ need not even be bounded, so the singularity of $Z_y(x)$ need not be $O(|x-y|^{2-n})$ as it was in the H\"older case.
These formulas simplify significantly if we use an affine change of variables in which $y$ corresponds to $x=0$ and $a_{ij}(0)=\de_{ij}$:
 \begin{equation}
Z(x)\sim |x|^{2-n}e^{I(|x|)} \quad\hbox{as}\ |x|\to 0,
\label{eq:Z-asym}
\end{equation}
where 
 \begin{equation}
I(r)=
\frac{1}{|S^{n-1}|}\int_{r<|z|<\e} \left({\rm tr}({\bf A}_z)-n\,\frac{\langle {\bf A}_z z, z\rangle}{|z|^2}\right)\frac{dz}{|z|^n}.
\label{def:I}
\end{equation}
We can verify that the absolute value of the integrand in (\ref{def:I}) is bounded by $\om(|z|)$, so the coefficients being Dini continuous implies that the improper integral defining $I(0)$ converges absolutely. Even if $I(0)$ does not converge, 
we shall see that for any $\lambda>0$ there exists $C_\lambda >0$ such that
 \begin{equation}
|I(r)|\leq \lambda |\log r|+C_\lambda \quad\hbox{for}\ 0<r<\e,
\label{eq:I(r)-est}
\end{equation} 
so the singularity of $Z$ at $x=0$ is never very far from $|x|^{2-n}$.
Nevertheless, the behavior of $I(r)$ as $r\to 0$ plays an important role for our results. There are three principal cases to consider.

\smallskip\noindent
{\it
1.  $I(0)=\lim_{r\to 0}I(r)$ exists and is finite.}

\noindent
In this case, we can scale $Z(x)$ by a constant multiple to obtain 
a solution  that is 
asymptotic to the fundamental solution for the Laplacian. In fact, we can solve the distributional equation
\begin{equation}
-{\mathcal L}(x,\del_x)Z(x)=C_0\,\delta(x),
\label{eq:-LZ=C_0delta}
\end{equation}
and find
\begin{equation}
C_0=(n-2)\,|S^{n-1}|\,e^{I(0)}.
\end{equation}
Note that the improper integral defining $I(0)$ may converge even if the modulus of continuity does {\it not} satisfy the Dini condition.

\eject
\smallskip\noindent
{\it
2.  $I(r)\to -\infty$ as $r\to 0$.}

\noindent
We see that $Z(x)= o(|x|^{2-n})$ as $|x|\to 0$;  we still have $Z(x)\to +\infty$ as $x\to 0$, so this violates the extended maximum principle of \cite{GS}. Nevertheless, we can solve (\ref{eq:-LZ=C_0delta}) to find  $C_0=0$ and this yields the surprising fact that there exists a distributional solution of ${\mathcal L}(x,\del_x)Z(x)=0$ for $x\in B_\e(y)$ which is not even bounded.

\smallskip\noindent
{\it
3. $I(r)\to \infty$ as $r\to 0$.}

\noindent
Now we find $Z(x)|x|^{n-2}\to \infty$ as $|x|\to 0$, so this solution grows more rapidly than the fundamental solution for the Laplacian. However, we cannot solve (\ref{eq:-LZ=C_0delta}) for $C_0$  in this case.

\smallskip
We next allow $y$ to vary over $U$. Provided that we are in Case 1 at each $y\in U$, we can use $Z_y(x)$ to construct the fundamental solution $F(x,y)$ in $U$. Our main result (Theorem 3) states that, provided $a_{ij}\in C^\om(U)$ where $\om$ satisfies (\ref{eq:om1}) and $I_y(0)=\lim_{r\to 0}I_y(r)$ exists and is finite for each $y\in U$, a fundamental solution $F(x,y)$ exists in the form (\ref{eq:F-asym-Holder}) where the remainder term $H(x,y)$ may be estimated in $L^p$ for any $p\in (1,\infty)$ in terms of $\om$ and the rate of convergence $I_y(r)\to I_y(0)$. More specifically, let us assume $\om(r)\,r^{-1+\kappa}$ is nonincreasing for $0<r<1$ where $\kappa\in (0,1)$, and introduce
\begin{equation}
 \sigma(r)=\int_0^r \om^2(t)\,\frac{dt}{t}.
\label{def:tilde-om}
\end{equation}
Further, let us  assume that
\begin{equation}
|I_y(r)-I_y(0)|\leq \theta(r) \quad \hbox{for all} \ y\in U,
\label{eq:theta}
\end{equation}
where $\theta(r)$ is a nondecreasing function of $r$ satisfying $\theta(0)=0$.
 Then for $y$ in a compact subset of $U$, we have
 \begin{equation}
 r^2M_p(D^2 H(\cdot,y),r;y)\leq c\,\max(\om(r),\sigma(r),\theta(r)) \quad\hbox{as}\ r\to 0,
 \label{est:Mp(D^2H)}
 \end{equation}
where $M_p(f(\cdot,y),r;y)$ denotes the $L^p$-mean of $f(x,y)$ as a function of $x$ (for fixed $y$) over the annulus $A_r(y)=\{x:r<|x-y|<2r\}$. 
(In fact, when $\om$ satisfies the Dini condition, the right hand side of (\ref{est:Mp(D^2H)}) reduces to just $c\,\om(r)$.)
Taking $p>n$, we obtain pointwise bounds on $|H(x,y)| + r|D_xH(x,y)|$, but we no longer have pointwise bounds on $D^2_xH(x,y)$ as we did in (\ref{eq:H-est}) when the coefficients $a_{ij}$ were H\"older continuous.

\smallskip\noindent
{\bf Note.} In (\ref{est:Mp(D^2H)}) and throughout this paper, $c$ is used to denote a constant whose value may change line-by-line. It may depend upon $n$ and the $a_{ij}$, but not on $r$.

\smallskip\noindent
{\bf Organization of this paper and comparison with other works.}
The 
organization of this paper is as follows. In Section 1 we discuss some preliminary estimates for solutions of the Poisson equation $\Lap u=f$. In Section 2 we construct the solution $Z_y(x)$ of (\ref{eq:LZ_y=0}) by first considering the case $a_{ij}(0)=\de_{ij}$ and then performing a change of coordinates. In Section 3 we analyze the equation $-{\mathcal L}(x,\del_x)Z_y(x)=C_y\delta(x-y)$ and calculate $C_y$ when $I_y(0)$ is finite or $-\infty$. Finally, in Section 4 we apply the previous results to construct our fundamental solution $F(x,y)$ in the form ({\ref{eq:F-asym-Holder}) and obtain estimates on the remainder term $H(x,y)$. The analysis in each section makes use of $L_p$-means on annuli $A_r=\{x:r<|x|<2r\}$ to measure the growth and decay of functions as $r\to 0$; these are defined in Section 1.

The results obtained and the techniques used in this paper are closely related to those
 in our previous paper ~\cite{MM}. In ~\cite{MM}  we studied the asymptotics for solutions of the 
 adjoint equation for ${\mathcal L}$, and discussed their relationship to the estimates obtained in ~\cite{E} as $y\to x$. In fact, the condition that $I_y(r)$ is bounded (above and below) as $r\to 0$ not only allows one to conclude that the solution $Z_y$ of (\ref{eq:LZ_y=0}) satisfies 
 $c_y|x-y|^{2-n}\leq Z_y(x)\leq c'_y|x-y|^{2-n}$ as $x\to y$, but that ${\mathcal L}^*(x,\del_x)u(x)=0$
 admits solutions in $B_\e(y)$ whose $L_p$-mean is bounded between positive constants as $r=|x-y|\to 0$; if $I_y(0)$ exists and is finite for every $y$, then this bounded solution of ${\mathcal L}^*(x,\del_x)u(x)=0$ is continuous.
 
Both ~\cite{MM} and the present work are independent of, but related to, the asymptotic theory developed in ~\cite{KM3}. In particular, $L_p$-means were extensively used in ~\cite{KM1} and 
~\cite{KM2}. The asymptotic formulas that we obtain are analogous to those of ~\cite{KM4}, where an asymptotic representation near the boundary was obtained for solutions to the Dirichlet problem for elliptic equations in divergence form with discontinuous coefficients. In particular, note that we do not make use of the maximum principle in this paper; in fact, most results of this paper hold for complex-valued coefficients, although this introduces some technical difficulties which we have chosen to avoid.


\section{Preliminaries}

Throughout this paper, to measure the growth of functions as $x\to y$, it will be convenient to use the $L^p$-mean for some $p\in (1,\infty)$:
 \begin{equation}
M_{p}(w,r;y)=\left(\meanint_{A_r(y)}|w(x)|^p\,dx\right)^{1/p},
\label{eq:Mpy}
\end{equation}
where $A_r(y)=\{x:r<|x-y|<2r\}$; here (and elsewhere in this paper) the slashed integral denotes mean value. It is natural to also introduce
 \begin{equation}
M_{\infty}(w,r;y)=\sup_{A_r(y)}|w(x)|.
\label{eq:M_infty}
\end{equation}
We may apply (\ref{eq:Mpy}) and (\ref{eq:M_infty}) to vector or matrix valued functions $w$ with $|w|$ denoting the norm. We also need to estimate derivatives, so we define
\begin{equation}
M_{1,\infty}(w,r;y)=rM_\infty(Dw,r;y)+M_\infty(w,r;y),
\label{eq:M1infty}
\end{equation}
where $Dw$ represents the gradient of $w$, and for $p\in(1,\infty)$
\begin{equation}
M_{2,p}(w,r;y)=r^2M_p(D^2w,r;y)+rM_{p}(Dw,r;y)+M_p(w,r;y),
\label{eq:M2p}
\end{equation}
where $D^2w$ represents the Hessian matrix of $w$.  Sobolev estimates show that
\begin{equation}
p>n \quad \Rightarrow \quad 
M_{1,\infty}(w,r;y)\leq c\, M_{2,p}(w,r;y).
\label{eq:Sobolev}
\end{equation}
When $y=0$, we shall abbreviate $M_p(w,r;0)$ as $M_p(w,r)$ (and similarly for $M_{1,\infty}$ and $M_{2,p}$).

For $x\in \RR^n\backslash\{0\}$, let $\theta=x/|x|\in S^{n-1}$ and let $d\theta$ denote the standard surface measure on $S^{n-1}$. 
We will use the spherical mean of a function $w$:
\begin{equation}
\overline{w}(r)=\meanint_{ S^{n-1}} w(r\theta)\,d\theta.
\label{eq:sph-mean}
\end{equation}
In particular, in this section we consider the equation
\begin{equation}
\Lap v=f \qquad\hbox{in}\ \RR^n\backslash\{0\}
\label{eq:Laplace}
\end{equation}
when $\overline f=0$, and investigate the behavior of the $L^p$-mean of the solution as $x\to 0$;   
our results are quite analogous to those of ~\cite{KM1} and ~\cite{KM3}. 
We shall let $\Gamma(|x|)=c_n\,|x|^{2-n}$ denote the fundamental solution for the Laplacian in $\RR^n$.

\begin{proposition}
Suppose $n\geq 2$, $p\in (1,\infty)$, and $f\in L_{\loc}^p(\RR^n\backslash\{0\})$ satisfies
$\overline f=0$, 
\begin{equation}
\int_{|x|< 1} |x|\,|f(x)|\,dx<\infty, 
\quad\hbox{and}\quad
\int_{|x|> 1} \frac{|f(x)|}{|x|^{n-1}}\,dx<\infty.
\label{eq:f_at_infty}
\end{equation}
Then $v=Kf=\Gamma\star f$ defines a distribution solution of (\ref{eq:Laplace}) that satisfies
\begin{equation}
M_{2,p}(Kf,r)\leq c
\left(r^2 \tilde M_p(f,r)+r^{1-n}\int_{|x|<r}|x|\,|f(x)|\,dx+r\int_{|x|>r}\frac{|f(x)|}{|x|^{n-1}}\,dx
\right),
\label{eq:M-est1-v}
\end{equation}
where we have introduced
\[
\tilde{M}_p(w,r):=\left(\meanint_{r/2<|x|<4r}|w(x)|^p\,dx\right)^{1/p}.
\]
\label{pr:1}
\end{proposition}

Elementary estimates can be used to show
\begin{equation}
\int_{|x|<r}|g(x)|\,dx\leq c\,\int_0^r M_p(g,\rho)\rho^{n-1}\,d\rho
\label{est:int-Mp}
\end{equation}
with an analogous estimate for $|x|>r$. These estimates enable us to obtain the following corollary from Proposition \ref{pr:1}.

\begin{cor}
Under the conditions of Proposition \ref{pr:1},
\begin{equation}
M_{2,p}(Kf,r)\leq c\left(r^{1-n}\int_0^r M_p(f,\rho)\rho^{n}\,d\rho
+r\int_r^\infty M_p(f,\rho)\,d\rho.\right)
\label{eq:M-est2-v}
\end{equation}
\end{cor}

\noindent
{\bf Proof of Proposition 1.} 
Let $v=Kf$ and let $\chi$ denote the characteristic function for the annulus $\tilde A_r=B_{4r}\backslash B_{r/2}$. Using $\overline{f}=0$, we see that
\begin{eqnarray}
v(x)=\int_{\RR^n}\Gamma(|x-y|)\chi(y)f(y)\,dy+
\int_{|y|<r/2}(\Gamma(|x-y|)-\Gamma(|x|))f(y)\,dy \nonumber \\
+\int_{|y|>4r}(\Gamma(|x-y|)-\Gamma(|y|))f(y)\,dy.\nonumber
\label{eq:v(x)=}
\end{eqnarray}
We want to estimate $M_p(v,r)$, $rM_p(Dv,r)$, and $r^2 M_p(D^2v,r)$, and show that each is bounded by the right hand side of (\ref{eq:M-est1-v}).

Using Stein's inequality \cite{St1}, for $0<a<n/p$ and $0<b<n/p'$ with $a+b=2$ we have
\[
\left\|\int_{\RR^n}\Gamma(|x-y|)\chi(y)f(y)\,dy\right\|_{L^p(A_r)}
\leq c\,r^a\|\chi(y)|y|^bf(y)\|_{L^p(\RR^n)}
\]
\[
= c\,r^a\||y|^b f(y)\|_{L^p(\tilde A_r)}\leq c\,r^2\|f\|_{L^p(\tilde A_r)}.
\]
It is elementary to show that for $|y|<r/2$ and $r<|x|<2r$ we have
$|\Gamma(|x-y|)-\Gamma(|x|)|\leq c\,|x|^{1-n}|y|$, where $c$ is independent of $r$, so
for $x\in A_r$ we have
\[
\left|
\int_{|y|<r/2}(\Gamma(|x-y|)-\Gamma(|x|))f(y)\,dy
\right|\leq 
c\,|x|^{1-n}\int_{|y|<r/2}|f(y)||y|\,dy.
\]
Similarly, we can show that for $|y|>4r$ and $r<|x|<2r$ we have
$|\Gamma(|x-y|)-\Gamma(|y|)|\leq c|x||y|^{1-n}$, so  $x\in A_r$  implies
\[
\left|
\int_{B_{4r}^c}(\Gamma(|x-y|)-\Gamma(|y|))f(y)\,dy\right|
\leq c\,|x|\int_{|y|>4r}\frac{|f(y)|}{|y|^{n-1}}\,dy.
\]
Using these estimates, we easily conclude that $M_p(v,r)$ is bounded by the right hand side of (\ref{eq:M-est1-v}).

Next we consider
\[
\del_i v(x) =\!\int_{\RR^n}\! \Gamma_i(x-y)\chi(y)f(y)\,dy
+\int_{B_{r/2}}\!\Gamma_i(x-y)f(y)\,dy+\int_{B^c_{4r}}\!\Gamma_i(x-y)f(y)\,dy,
\]
where
\[
\Gamma_i(x)=\Gamma'(|x|)\frac{x_i}{|x|}.
\]
Applying Stein's inequality as above but with $a+b=1$, we conclude
\[
\left\|\int_{\RR^n} \Gamma_i(x-y)\chi(y)f(y)\,dy\right\|_{L^p( A_r)}
\leq c\,r\,\|f\|_{L^p(\tilde A_r)}.
\]
Elementary estimates for $r<|x|<2r$ show that
\[
\left|\int_{B_{r/2}}\!\Gamma_i(x-y)f(y)\,dy\right| =\left|\int_{B_{r/2}}\!(\Gamma_i(x-y)-\Gamma_i(x))f(y)\,dy\right| 
\leq c\,|x|^{-n}\int_{|y|<r/2}|f(y)||y|\,dy
\]
and
\[
\left|
\int_{B_{4r}^c}\!\Gamma_i(x-y)f(y)\,dy\right|
=\left|
\int_{B_{4r}^c}\!(\Gamma_i(x-y)-\Gamma_i(y))f(y)\,dy\right|
\leq c\,\int_{|y|>4r}\!\frac{|f(y)|}{|y|^{n-1}}\,dy.
\]
From these estimates we easily conclude that $rM_p(Dv,r)$ is bounded by the right hand side of 
(\ref{eq:M-est1-v}).

Finally, we consider
\[
\del_i\del_j v(x)=
\int_{\RR^n}\!\Gamma_{ij}(x-y)\chi(y)f(y)\,dy
+\int_{B_{r/2}}\!\Gamma_{ij}(x-y)f(y)\,dy
+\int_{B_{4r}^c}\!\Gamma_{ij}(x-y)f(y)\,dy,
\]
where $\Gamma_{ij}$ is the singular kernel given by
\[
\Gamma_{ij}(x)=\Gamma''(|x|)\frac{x_ix_j}{|x|^2}+\Gamma'(|x|)\frac{\de_{ij}|x|^2-x_ix_j}{|x|^3}.
\]
 Using the $L^p$-boundedness of singular integral operators, we conclude
\[
\left\|\int_{\RR^n}\Gamma_{ij}(x-y)\chi(y)f(y)\,dy\right\|_{L^p(A_r)}
\leq c\,\|f\|_{L^p(\tilde A_r)}.
\]
Elementary estimates for $r<|x|<2r$ show that
\[
\left|\int_{B_{r/2}}\Gamma_{ij}(x-y)f(y)\,dy\right|=\left|\int_{B_{r/2}}(\Gamma_{ij}(x-y)-\Gamma_{ij}(x))f(y)\,dy\right|
\]
\[
 \leq c\,|x|^{-n-1}\int_{|y|<r/2}|f(y)||y|\,dy
\]
and
\[
\left|
\int_{B_{4r}^c}\Gamma_{ij}(x-y)f(y)\,dy\right|
=\left|
\int_{B_{4r}^c}(\Gamma_{ij}(x-y)-\Gamma_{ij}(|y|))f(y)\,dy\right|
\]
\[
\leq c\,|x|^{-1}\int_{|y|>4r}\frac{|f(y)|}{|y|^{n-1}}\,dy.
\]
These estimates show that $r^2M_p(D^2v,r)$ is bounded by the right hand side of (\ref{eq:M-est1-v}).
This completes the proof.
\hfill $\Box$

\section{Constructing the Singular Solution $Z_y(x)$ in $B_\e(y)\backslash\{y\}$}

In this section, we fix $y$ and construct a solution of (\ref{eq:LZ_y=0}) for $\e$ sufficiently small. 
Since the result is localized near $y$, the ellipticity and continuity of the coefficients of ${\mathcal L}$ need only be assumed at $y$. In fact, 
we first consider the case when $y=0$ and $a_{ij}(0)=\de_{ij}$: 
\begin{equation}
\sup_{ |x|=r}\|{\bf A}_x-{\bf I}\|\leq
\om(r)\quad\hbox{for}\ 0<r<1,
\label{eq:aij(0)-asym}
\end{equation}
 where ${\bf I}$ is the identity matrix and $\om$ is a continuous, nondecreasing function for $0\leq r< 1$ satisfying the square-Dini condition (\ref{eq:om1}). 
 For convenience, we shall also assume the monotonicity condition that for some $\kappa\in(0,1)$ we have
\begin{equation}
\om(r)\,r^{-1+\kappa}\quad\hbox{is nonincreasing for $0<r<1$}.
\label{eq:om2}
\end{equation}
The significance of (\ref{eq:om2}) is that it requires $\om(r)$ to vanish more slowly than $r$ as $r\to 0$, which is a natural constraint. (As long as $\om(r)\,r^{-1+\kappa}$ is nonincreasing for $0<r<\e$ with some $\e>0$, then $\om$ may be modified for $\e<r<1$ to satisfy (\ref{eq:om2}).) We seek a solution of
\begin{equation}
{\mathcal L}(x,\del_x)\,Z(x)=0 \quad\hbox{for}\ x\in B_\e\backslash\{0\},
\label{eq:LZ=0}
\end{equation}
where $B_\e=\{x:|x|<\e\}$ with $\e$  sufficiently small. Although we generally assume that ${\mathcal L}$ has real-valued coefficients, the theorem below holds when the $a_{ij}$ are complex-valued.

\begin{theorem}
For $n\geq 3$ and $p\in (1,\infty)$, suppose the 
coefficients $a_{ij}(x)$ are bounded, measurable functions satisfying (\ref{eq:aij(0)-asym}). For  $\e>0$ sufficiently small, there exists a solution  of (\ref{eq:LZ=0}) in the form 
 \begin{equation}
  Z(x)=h(|x|)+v(x),
  \label{eq:Z}
  \end{equation}
    where $h$ is of the form
     \begin{equation}
 h(r)=\int_r^\e s^{1-n} \,e^{I(s)}\,ds\left(1+\zeta(r)\right),
  \label{est:h}
 \end{equation}
 with $I(r)$ given by (\ref{def:I}) and
 \begin{equation}
M_{2,p}(\zeta,r)\leq c\,\max(\om(r),\sigma(r)),
\label{eq:M2p(zeta)}
\end{equation}
where $ \sigma$ is given in (\ref{def:tilde-om}), and $v$ in (\ref{eq:Z}) satisfies
\begin{equation}
M_{2,p}(v,r)\leq c\,r^{2-n}\,e^{I(r)}\,\om(r).
\label{eq:M2p(v)}
\end{equation}
Moreover, for  any  $u\in W_{\loc}^{2,p}(\overline{B_{\e}}\backslash\{0\})$ that is a strong
solution of ${\mathcal L}(x,\del_x)u=0$ in $\overline{B_{\e}}\backslash\{0\}$ subject to the growth condition
\begin{equation}
M_{2,p}(u,r)\leq c\, r^{1-n+\e_0} \quad\hbox{where} \ \e_0>0,
\label{eq:u-hyp}
\end{equation}
 there exist constants $C, C_0, C_1,\dots,C_n$ (depending on $u$) such that
\begin{equation}
u(x)=CZ(x)+C_0+\sum_{j=1}^n C_jx_j+w(x),
\label{eq:u-asym}
\end{equation}
where $w$ satisfies
\begin{equation}
M_{2,p}(w,r)\leq c\,r^{2-\e_1} \quad\hbox{for any}\ \e_1>0.
\label{eq:M2p(w)}
\end{equation}
\label{th:I}
\end{theorem}

We shall prove this theorem below, but first let us make some observations. 
In general, we do not know whether $I(r)$ is bounded as $r\to 0$, but we can verify that
$|I'(r)|\leq c\,r^{-1}\om(r)$, so
integration by parts in (\ref{est:h}) shows that
\begin{equation}
h(r)=\frac{r^{2-n}}{n-2}\,e^{I(r)} + h_1(r), 
\label{eq:h-asym}
\end{equation}
where $h_1(r)$ satisfies 
$M_{1,\infty}(h_1,r)\leq c\,r^{2-n}e^{I(r)}\max(\om(r),\sigma(r)).$ 
If we take $p>n$ and apply (\ref{eq:Sobolev}) to $v$, we  conclude that
\begin{equation}
Z(x)=\frac{|x|^{2-n}e^{I(|x|)}}{n-2}\left(1+\xi(x)\right) \quad\hbox{as}\ |x|\to 0.
\label{eq:Z-asym2}
\end{equation}
where $M_{1,\infty}(\xi,r)\leq c\, \max(\om(r),\sigma(r)).$ 
Obviously, we can multiply the $Z$ of (\ref{eq:Z-asym2}) by $n-2$ to obtain the 
$Z$ of (\ref{eq:Z-asym}). 

Even when $I(r)$ is not bounded as $r\to 0$, we can derive useful bounds on $Z(x)$ as $|x|\to 0$. 
It is not difficult to verify that a symmetric matrix ${\bf A}$ satisfies
\begin{equation}
-2(n-1)\|{\bf A}-{\bf I}\|\leq {\rm tr}({\bf A})-n\langle {\bf A}y,y\rangle\,|y|^{-2}
\leq 2(n-1)\|{\bf A}-{\bf I}\| \quad\hbox{for}\ |y|=1,
\label{eq:integrand-bound}
\end{equation}
so there exist constants $c, C>0$ so that $Z$  satisfies (assuming $n\geq 3$)
 \begin{equation}
c\,|x|^{2-n}\,\exp\left(-c_n\int_{|x|}^\e\om(t)\frac{dt}{t}\right)
\leq |Z(x)|
\leq C\,|x|^{2-n}\,\exp\left(c_n\int_{|x|}^\e\om(t)\frac{dt}{t}\right)
\label{eq:Z-est1}
\end{equation} 
as $x\to 0$, where $c_n=2(n-1)/|S^{n-1}|$.
Using (\ref{eq:Z-est1}) and the fact that $\om(r)\to 0$ as $r\to 0$, we obtain (\ref{eq:I(r)-est}),
which shows that the singularity of $Z$ at $x=0$ is very close to the classical case.

An interesting class of examples is obtained by letting
\begin{equation}
a_{ij}(x)=\de_{ij}+g(|x|)x_ix_j|x|^{-2},
\label{eq:GS1}
\end{equation}
where $g(0)=0$ but vanishes slowly as $r\to 0$. Gilbarg \& Serrin \cite{GS} used (\ref{eq:GS1}) with certain specific functions $g$  to show that Dini continuity is essential for their extended maximum principle to hold. In our formulation,
\begin{equation}
{\rm tr}({\bf A}_z)-n\,\frac{\langle {\bf A}_z z, z\rangle}{|z|^2}=(1-n)g(|z|),
\label{eq:GS2}
\end{equation}
so
\begin{equation}
I(r)=(1-n)\int_r^\e g(\rho)\,\frac{d\rho}{\rho}.
\label{eq:GS3}
\end{equation}
Thus any $g(r)>0$ which does {\it not} satisfy the Dini condition at $r=0$ (but does satisfy the square-Dini condition) will yield $I(r)\to -\infty$ as $r\to 0$, so the $Z(x)$ of Theorem 1 gives an example of a solution of (\ref{eq:LZ=0}) with singularity at $x=0$ even though $Z(x)=o(|x|^{2-n})$ as $|x|\to 0$, i.e.\  the extended maximum principle fails; the specific function in \cite{GS} is $g(r)=-(1+ (n-1)\log r)^{-1}$.

\smallskip
\noindent{\bf Proof of Theorem 1.}
Instead of showing the existence of $Z$ in a very small ball $B_\e$, we shall replace the condition that $\om$ satisfies the square-Dini condition by
\begin{equation}
\sigma(1)=\int_0^1 \om^2(t)\,\frac{dt }{t}< \de,
\label{eq:om3}
\end{equation}
where $\de$ is sufficiently small, and show existence in the unit ball $B_1$.
In fact, using (\ref{eq:om3}) and  (\ref{eq:om2}), we see that
\begin{equation}
\omega(r)< c_{\kappa,n}\sqrt{\de} \quad\hbox{for}\ 0<r<1,
\label{eq:om^2<de}
\end{equation}
where $c_{\kappa,n}$ depends only on $\kappa$ and $n$:
\[
\de>\int_{r/2}^r \om^2(t)\,\frac{dt }{t}\geq \om^2(r)\,r^{-2+2\kappa}\int_{r/2}^r t^{1-2\kappa}\,dt= \om^2(r)\ c'_{\kappa,n}.
\]
Moreover,  it will be useful to consider ${\mathcal L}$ as defined on all of $\RR^n$ with
${\mathcal L}=\Delta$ outside of $B_1$. Therefore, we shall assume that
\begin{equation}
a_{ij}(x)=\de_{ij} \qquad\hbox{for} \  |x|>1,
\label{eq:aij=dij}
\end{equation} 
and investigate a solution of ${\mathcal L}Z=0$ in $\RR^n\backslash\{0\}$.

To construct $Z(x)$, we let $h(r)=\overline{Z}(r)$ denote the spherical mean as in (\ref{eq:sph-mean}), and let  $v(x)=Z(x)-h(|x|)$, so $\overline{v}(r)=0$. We shall reduce the problem to solving an operator equation of the form $(I+S+T)v=f$ where $S$ and $T$ have small operator norm on a Banach space $X$ defined as follows: for fixed $p\in (1,\infty)$, let us  consider the functions $v$ in $W^{2,p}_{\loc}(\RR^n\backslash\{0\})$ for which the norm
\begin{equation}
\|v\|_{X}=\sup_{0<r<1}\frac{M_{2,p}(v,r)r^{n-2}}{\om(r)e^{I(r)}}+\sup_{r>1}\frac{M_{2,p}(v,r)\,r^{n-1}}{\sqrt{\de}}
\label{def:X}
\end{equation}
is finite, and take the closure to form a Banach space $X$. The reduction to the operator equation will take several steps.

If we take the spherical mean of ${\mathcal L}(h+v)=0$, we obtain an ordinary differential equation for $h$:
\[
\alpha(r)h''+\frac{\alpha_n(r)-\alpha(r)}{r}h'+\overline{a_{ij}\del_i\del_j v}(r)=0,
\]
where
\begin{equation}
\alpha_n(r)=\meanint_{S^{n-1}} a_{ii}(r\theta)\, d\theta \quad\hbox{ and}\quad
\alpha(r)=\meanint_{S^{n-1}} a_{ij}(r\theta)\theta_i\theta_j\,d\theta.
\label{def:alpha}
\end{equation}
From (\ref{eq:aij(0)-asym}) it is clear that these functions satisfy
\begin{equation}
|\alpha_n(r)-n|\leq \om(r)\ \hbox{and}\ 
|\alpha(r)-1|\leq \om(r)\ \hbox{for}\ 0<r<1.
\label{eq:alpha-asym}
\end{equation}
Using the fact that $\overline{\Lap  v}=\Lap \overline v=0$, the equation for $h$ becomes
\begin{equation}
\alpha(r)h''+\frac{\alpha_n(r)-\alpha(r)}{r}h'+\overline{\beta_{ij}\del_i\del_j v}(r)=0,
\label{eq:ode}
\end{equation}
where $\beta_{ij}(x)=a_{ij}(x)-\de_{ij}$ satisfies $|\beta_{ij}(x)|\leq \omega(|x|)$. Notice that $v$ satisfies the partial differential equation
\begin{equation}
-\Lap v=\beta_{ij}\del_i\del_j h-\overline{\beta_{ij}\del_i\del_j h} +
\beta_{ij}\del_i\del_j v -\overline{\beta_{ij}\del_i\del_j v}.
\label{eq:pde-v}
\end{equation}
We shall first solve (\ref{eq:ode}) and use that to eliminate $h$ from (\ref{eq:pde-v}); then we will be able to apply the results of Section 1. 

Let us introduce $g=h'$ and rewrite (\ref{eq:ode}) as
\begin{equation}
g'+\frac{n-1+R(r)}{r}g=B[D^2v](r),
\label{eq:g}
\end{equation}
where
\begin{equation}
R(r)=\frac{\alpha_n(r)}{\alpha(r)}-n
\label{def:R}
\end{equation}
and
\[
B[D^2 v](r)=-\alpha^{-1}(r)\overline{\beta_{ij}\del_i\del_j v}(r).
\]
Notice that
\begin{equation}
|R(r)|\leq c\,\omega(r)\ \hbox{for}\ 0<r<1 \quad\hbox{and}\quad R(r)=0 \ \hbox{for}\ r>1,
\label{eq:R-om}
\end{equation}
and
\[
|B[D^2 v](r)|\leq c\,\omega(r)|\overline{D^2v}(r)|\ \hbox{for}\ 0<r<1 
 \quad\hbox{and}\quad B[D^2 v](r)=0 \ \hbox{for}\ r>1.
\]
Moreover, the monotonicity of $\om(r)$ together with (\ref{eq:om2}) imply
\begin{equation}
\max_{r\leq\rho\leq 2r}\om(\rho)\leq c\,\om(r),
\label{eq:om(2r)}
\end{equation}
so we consequently  obtain
\begin{equation}
M_p(B[D^2v],r)\leq c\,\omega(r)M_p(D^2v,r)\quad\hbox{for}\ 0<r<1.
\label{eq:Mp-B}
\end{equation}

To solve (\ref{eq:g}), let us introduce 
\begin{equation}
E_\pm(r)=\exp\left[\pm\int_r^\infty R(t)\,\frac{dt}{t}\right]=\exp\left[\pm\int_r^1 R(t)\,\frac{dt}{t}\right]=\frac{1}{E_\mp(r)}
\label{def:E2}
\end{equation}
and observe that $E_\pm(r)\equiv 1$ for $r>1$. It is useful to 
observe that
\[
E_-(r)E_+(\rho)=\exp\left(\int_\rho^r R(t)\,\frac{dt}{t}\right),
\]
so as a consequence of (\ref{eq:R-om}) and (\ref{eq:om^2<de}), we obtain
\begin{equation}
\left(\frac{\rho}{r}\right)^{c\sqrt{\de}}\leq \exp\left(\pm\int_\rho^r R(t)\,\frac{dt}{t}\right)\leq \left(\frac{r}{\rho}\right)^{c\sqrt{\de}}
\quad\hbox{for}\ 0<\rho\leq r\leq 1.
\label{eq:E-est}
\end{equation}
 In particular, 
 \begin{equation}
 c_1\,E_\pm(r)\leq E_\pm(\rho)\leq c_2\,E_\pm(r) \quad\hbox{for}\ r<\rho<2r,
 \label{eq:E(2r)}
 \end{equation}
 and 
 for any $g\in  L^p_\loc (\RR^n\backslash\{0\})$ we can readily verify that
 \begin{equation}
 M_p(|x|^\nu\,E_{\pm}(|x|)\,g(x),r)\leq c\,r^\nu\,E_\pm(r)\,M_p(g,r),
 \label{est:Mp(Eg)}
 \end{equation}
for any fixed $\nu\in \RR$.

It will be more convenient for us to use $E_{\pm}(r)$ than $e^{\pm I(r)}$, but these functions are equivalent: if we note that (\ref{def:I}) can be written as 
\begin{equation}
I(r)=\int_r^1[\alpha_n(\rho)-n\alpha(\rho)]\,\frac{d\rho}{\rho},
\label{eq:I}
\end{equation}
then we see that
\begin{equation}
E_+(r)=A\,e^{I(r)}(1+\tau(r)),
\label{def:A}
\end{equation}
where $A=\exp[\int_0^1R(\rho)[1-\alpha(\rho)]\rho^{-1}{d\rho}]$ is finite and positive,  
and 
\begin{equation}
\tau(r)=\exp\left[-\int_0^r R(\rho)(1-\alpha(\rho))\,\frac{d\rho}{\rho}\right]-1
\label{def:tau}
\end{equation}
satisfies $|\tau(r)|\leq c\,\sigma(r)$. Thus for some constants $c_1,c_2$ we have
\begin{equation}
c_1 E_+(r)\leq e^{I(r)}\leq c_2 E_+(r)\quad\hbox{for}\ 0<r<1.
\label{est:E+e^I}
\end{equation}

Now if we introduce
$\phi(r)=r^{n-1} E_-(r) g(r),$
then we can rewrite (\ref{eq:g}) as
\begin{equation}
\phi'(r)=r^{n-1} E_-(r) B[D^2v](r).
\label{eq:phi'}
\end{equation}
But (\ref{eq:phi'}) may be integrated to obtain
\begin{equation}
\phi(r)=\phi(0)+\int_0^r\rho^{n-1} E_-(\rho) B[D^2v](\rho)\,d\rho,
\label{eq:phi}
\end{equation}
where $\phi(0)$ is an arbitrary constant. Of course, to conclude (\ref{eq:phi}), we must verify that $\phi'$ is integrable on $(0,1)$. But $v\in X$ implies $M_p(D^2 v,r)\leq c\,\omega(r)r^{-n}E_+(r)$,
so we can use (\ref{eq:om(2r)}), (\ref{eq:Mp-B}), (\ref{eq:E(2r)}),  and H\"older's inequality to obtain
\[
\left|\int_r^{2r}\rho^{n-1}E_-(\rho) B[D^2v](\rho)\,d\rho\right|
\leq c\,E_-(r)\om(r)\int_{r<|x|<2r}|D^2v(x)|dx
\leq c\, \om^2(r).
\]
Now if we write 
\[
\int_0^s \rho^{n-1}E_-(\rho)B[D^2v](\rho)\,d\rho=
\sum_{j=0}^\infty \int_{s/2^{j+1}}^{s/2^j} \rho^{n-1}E_-(\rho)B[D^2v](\rho)\,d\rho,
\]
then we obtain the estimate 
\begin{equation}
\begin{aligned}
\left| \int_0^s \rho^{n-1}E_-(\rho)B[D^2v](\rho)\,d\rho \right|
&\leq c\sum_{j=0}^\infty \om^2\left(\frac{s}{2^{j+1}}\right)  \\
&\leq c\int_0^s \om^2(\rho)\,\frac{d\rho}{\rho}=c\,\sigma(s)<c\,\de. 
\end{aligned}
\label{eq:int-phi'}
\end{equation}
In particular, we see that $\phi'$ is integrable and so (\ref{eq:phi}) is valid.
This enables us to write
\begin{equation}
h'(r)=g(r)=r^{1-n}E_+(r)\left[\phi(0)+\int_0^r\rho^{n-1} E_-(\rho) B[D^2v](\rho)\,d\rho\right]
\label{eq:h'}
\end{equation}
and
\begin{equation}
h''(r)=\frac{1-n-R(r)}{r^{n}} E_+(r)\left[\phi(0)+\int_0^r\rho^{n-1}E_-(\rho)B[D^2v](\rho)\,d\rho\right]+B[D^2](r).
\label{eq:h''}
\end{equation}

We can now use (\ref{eq:h'}) and (\ref{eq:h''}) to
express
\[
\beta_{ij}\del_i\del_j h=r^{-n}E_+(r)\!
\left[\phi(0)+\int_0^r\!\rho^{n-1}E_-(\rho)B[D^2v](\rho)d\rho\right]\!\psi(r\theta)
\]
\[
+B[D^2v](r)\,\beta_{ij}\theta_i\theta_j,
\]
where 
\begin{equation}
\psi(r\theta)=\beta_{ii}(r\theta)-(n+R(r))\beta_{ij}(r\theta)\theta_i\theta_j
\label{def:psi}
\end{equation}
 also satisfies
$|\psi(r\theta)|\leq c\,\om(r)$.
 Thus we can apply $K=\Gamma\star$ to 
(\ref{eq:pde-v}) to obtain
\begin{equation}
v+Sv+Tv=\phi(0)w,
\label{eq:v-int}
\end{equation}
where
\[
w(x)=K_{y\to x}\,[|y|^{-n}E_+(|y|)(\psi(y)-\overline{\psi}(|y|)],
\]
\[
Sv(x)=-K_{y\to x}\left[|y|^{-n}E_+(|y|)\int_0^{|y|}\rho^{n-1}E_-(\rho) B[D^2v](\rho)\,d\rho\,(\psi(y)-\overline{\psi}(|y|)\right],
\]
and
\[
Tv=-K\,\left[B[D^2v]\,(\beta_{ij}\theta_i\theta_j-\overline{\beta_{ij}\theta_i\theta_j})
+\beta_{ij}\del_i\del_j v-\overline{\beta_{ij}\del_i\del_j v}\right].
\]
For a given value of $\phi(0)$, we can solve (\ref{eq:v-int}) provided we show that $w\in X$ and
both $S:X\to X$ and $T:X\to X$ have small operator norms.

To show $w\in X$, we must estimate $M_{2,p}(w,r)$ for $0<r<1$ and for $r>1$. For $0<r<1$ apply
Proposition 2 to $f(x)=|x|^{-n}E_+(|x|)(\psi(x)-\overline{\psi}(|x|))$ (which vanishes for $|x|>1$) to conclude
\[
M_{2,p}(w,r)\leq c\left(r^{1-n}\int_0^r E_+(\rho)\om(\rho)\,d\rho+r\int_r^1 E_+(\rho)\om(\rho)\rho^{-n}\,d\rho\right).
\]
Using (\ref{eq:E-est}) and the fact that $r^{1-\kappa}\om(r)$ is nondecreasing (since both $r^{1-\kappa}$ and $\om$ are), we find
\[
\int_0^r E_+(\rho)\om(\rho)\,d\rho\leq  E_+(r)\,r^{c\sqrt{\de}}\int_0^r\rho^{-c\sqrt{\de}}\om(\rho)\,d\rho
\]
\[
\leq E_+(r)\,\om(r)\,r^{c\sqrt{\de}+1-\kappa}\int_0^r \rho^{-c\sqrt{\de}-1+\kappa}\,d\rho.
\]
Taking $\de$ small enough that $\kappa-\sqrt{c\de}>0$, we obtain
\begin{equation}
\int_0^r E_+(\rho)\om(\rho)\,d\rho\leq c\,r\,E_+(r)\,\om(r).
\label{est:int_0^r}
\end{equation}
Using (\ref{eq:E-est}) and the fact that $r^{-1+\kappa}\,\om(r)$ is nonincreasing (by (\ref{eq:om2})), we find
\[
\int_r^1 E_+(\rho)\om(\rho)\rho^{-n}\,d\rho\leq E_+(r)\,r^{-c\sqrt{\de}}\int_r^1\rho^{c\sqrt{\de}-n}\om(\rho)\,d\rho 
\]
\[
\leq E_+(r)\,\om(r)\,r^{-c\sqrt{\de}-1+\kappa}\int_r^1 \rho^{c\sqrt{\de}+1-\kappa-n}\,d\rho.
\]
For $\de$ small enough that $n-2+\kappa-c\sqrt{\de}>0$, we have
\begin{equation}
\int_r^1 E_+(\rho)\om(\rho)\,\rho^{-n}\,d\rho\leq c\,r^{1-n}\,E_+(r)\,\om(r).
\label{est:int_r^1}
\end{equation}
Using (\ref{est:int_0^r}) and (\ref{est:int_r^1}), we obtain
\[
\frac{M_{2,p}(w,r)r^{n-2}}{\om(r)E_+(r)}
\leq  c \quad\hbox{for all }\ 0<r<1.
\]
We can then use (\ref{est:E+e^I}) to replace $E_+(r)$ by $e^{I(r)}$ as required in the norm for $X$ in (\ref{def:X}).
Meanwhile, for $r>1$ we use $\om(1)\leq c\,\sqrt{\de}$ and $E_+(\rho)\leq  \rho^{-c\sqrt{\de}}$ for $0<\rho<1$ to find
\[
M_{2,p}(w,r)\leq c\, r^{1-n}\int_0^1 E_+(\rho)\om(\rho)\,d\rho  \leq c\,\sqrt{\de}\, r^{1-n}\int_0^1 \rho^{-c\sqrt{\de}}\,d\rho
\leq c\,\sqrt{\de}\, r^{1-n},
\]
provided $\de$ is sufficiently small. Consequently,
\[
\frac{M_{2,p}(w,r)\,r^{n-1}}{\sqrt{\de}}\leq c \quad\hbox{for all}\ r>1,
\]
and this confirms that $w\in X$.

Next let us show that $S$ maps $X$ to itself with small operator norm. We suppose that $\|v\|_X\leq 1$ and estimate $M_{2,p}(Sv,r)$ separately for $0<r<1$ and for $r>1$. For $0<r<1$, we examine the proof
of (\ref{eq:int-phi'}) and observe that the condition $\|v\|_X\leq 1$ enables us to choose the constant $c$ independent of $v$. Thus the function
\[
f_1(y)=|y|^{-n}E_+(|y|)\int_0^{|y|} \rho^{n-1}E_-(\rho)\,B[D^2v](\rho)\,d\rho\,(\psi(y)-\overline{\psi}(|y|))
\]
 satisfies 
 \[
 M_p(f_1,r)\leq c\,\de\,E_+(r)\,\om(r)\,r^{-n}\ \hbox{for}\ 0<r<1
 \]
  and $M_p(f_1,r)=0$ for $r>1$. For $0<r<1$, we apply Proposition 2 to $Sv=-Kf_1$ to obtain 
\[
M_{2,p}(Sv,r)\leq c\,\de\left(r^{1-n}\int_0^r E_+(\rho)\,\om(\rho)\,d\rho+
r\int_r^1E_+(\rho)\,\om(\rho)\,\rho^{-n}\,d\rho\right),
\]
and then use (\ref{est:int_0^r}) and (\ref{est:int_r^1}) to conclude (for $\de$ sufficiently small)
\[
\frac{M_{2,p}(Sv,r)r^{n-2}}{\om(r)E_+(r)}\leq c\,\de\quad\hbox{for all }\ 0<r<1.
\]
On the other hand, for $r>1$, Proposition 2 implies  (for $\de$ sufficiently small)
\[
M_{2,p}(Sv,r)\leq c\,\de\, r^{1-n}\int_0^1E_+(\rho)\,\om(\rho)\,d\rho
\leq c\,\de^{3/2}\,r^{1-n}\int_0^1\rho^{-c\sqrt{\de}}\,d\rho\leq c\,\de^{3/2}\,r^{1-n}.
\]
Thus we have
\[
\frac{M_{2,p}(Sv,r)r^{n-1}}{\sqrt{\de}}\leq c\,\de \quad\hbox{for all }\ r>1.
\]
Combining these inequalities, we see that $S:X\to X$ has small operator norm.

Finally, we show that $T$ maps $X$ to itself with small operator norm. We suppose that $\|v\|_X\leq 1$ and estimate $M_{2,p}(Tv,r)$ separately for $0<r<1$ and for $r>1$. Notice that the function
\[
f_2=B[D^2v]\,(\beta_{ij}\theta_i\theta_j-\overline{\beta_{ij}\theta_i\theta_j})
\]
 satisfies 
 \[
 M_p(f_2,r)\leq \om(r) M_p(B[D^2v],r)\leq c\,\om^3(r)\,E_+(r)\,r^{-n} \  \hbox{for}\ 0<r<1,
 \]
 where $c$ is independent of $v$, and $M_p(f_2,r)=0$ for $r>1$. 
Similarly, the function
\[
f_3=\beta_{ij}\del_i\del_j v-\overline{\beta_{ij}\del_i\del_j v}
\]
 satisfies 
 \[
 M_p(f_3,r)\leq \om(r) M_p(D^2v,r)\leq \om^2(r)\,E_+(r)\,r^{-n} \  \hbox{for}\ 0<r<1,
 \]
  and $M_p(f_3,r)=0$ for $r>1$.
For $0<r<1$, we apply Proposition 2 to $Tv=-K(f_2+f_3)$ to obtain
\[
M_{2,p}(Tv,r)\leq c\left(
r^{1-n}\int_0^r\om^2(\rho)\,E_+(\rho)\,d\rho+r\int_r^1\om^2(\rho)\,E_+(\rho)\rho^{-n}\,d\rho\right).
\]
Using (\ref{eq:om^2<de}), (\ref{est:int_0^r}), and (\ref{est:int_r^1}),
\[
\frac{M_{2,p}(Tv,r)\,r^{n-2}}{\om(r)E_+(r)}\leq c\sqrt{\de}\quad\hbox{for all}\ 0<r<1.
\]
On the other hand, for $r>1$, we use 
(\ref{eq:om^2<de}) and (\ref{eq:E-est}) to estimate
\[
M_{2,p}(Tv,r)\leq c\,r^{1-n}\int_0^1\om^2(\rho)\,E_+(\rho)\,d\rho
\leq c\,\de \, r^{1-n}\int_0^1 \rho^{-c\sqrt{\de}}\,d\rho
\leq c\,\de\,r^{1-n}.
\]
Consequently,
\[
\frac{M_{2,p}(Tv,r)\,r^{n-1}}{\sqrt{\de}}\leq c\,\sqrt{\de}\quad\hbox{for all}\ r>1.
\]
Combining these estimates, we see that $T:X\to X$ has small operator norm.

Since both $S$ and $T$ have small operator norms on $X$, we conclude that 
(\ref{eq:v-int}) has a unique solution $v$, depending on the choice of the constant $c_1=\phi(0)$.
But once $c_1$ and $v$ are known, we obtain $g(r)$ from (\ref{eq:h'}),
and $h(r)$ by integration of $g(r)$:
\begin{equation}
h(r)=\int_r^\infty s^{1-n}E_+(s)\left[c_1+\int_0^s\rho^{n-1}E_-(\rho)B[D^2v](\rho)\,d\rho\right]ds+c_2.
\label{eq:h}
\end{equation}
where $c_2$ is an arbitrary constant. 
To obtain the desired solution of Theorem 1, we choose $c_1$ to enable us to replace $E_+(r)$ by $e^{I(r)}$ for $0<r<1$. Using (\ref{def:A}) we see that we should choose $c_1=A^{-1}$ and write $h(r)=h_0(r)+h_1(r)+c$ where
\begin{equation}
h_0(r)=\int_r^1 s^{1-n}e^{I(s)}\,ds
\label{def:h0}
\end{equation}
and (recalling $\tau$ from (\ref{def:tau}))
\begin{equation}
h_1(r)=\int_r^1 s^{1-n}e^{I(s)}\tau(s)ds+
\int_r^1s^{1-n}E_+(s)\int_0^s \rho^{n-1}E_-(\rho)B[D^2v](\rho)d\rho ds.
\label{def:h1}
\end{equation}
Now  integrate by parts to obtain
\[
h_0(r)=\frac{r^{2-n}}{n-2}\,e^{I(r)}+c+\frac{1}{n-2}\int_r^1 s^{2-n}e^{I(s)}I'(s)\,ds.
\]
But $|I'(s)|\leq c\,\om(s)/s$ and, similar to (\ref{eq:E-est}), we can show that
\[
e^{I(s)}e^{-I(r)}\leq \left(\frac{s}{r}\right)^{c\sqrt{\de}}\quad\hbox{for}\ s>r,
\]
so we may use (\ref{eq:om2}) to obtain
\[
\left|\int_r^1 s^{2-n}e^{I(s)}I'(s)\,ds\right|\leq
c\,\int_r^1s^{1-n}e^{I(s)}\om(s)\,ds\leq c\,r^{-c\sqrt{\de}}e^{I(r)}\int_r^1 s^{1-n+c\sqrt{\de}}\om(s)\,ds
\]
\[
\leq c\,\om(r)e^{I(r)}r^{-c\sqrt{\de}-1+\kappa}[r^{3-n+c\sqrt{\de}-\kappa}+1]\leq c\,r^{2-n}\om(r)e^{I(r)}
\]
provided $\de$ is sufficiently small. Thus we find
\begin{equation}
\left| h_0(r)-\frac{r^{2-n}}{n-2}e^{I(r)}\right| \leq c\,r^{2-n}e^{I(r)}\om(r)\qquad\hbox{for}\ 0<r<1.
\label{eq:asym-h0}
\end{equation}
To estimate $h_1$ we use $|\tau(s)|\leq c\,\sigma(s)$  together with (\ref{eq:int-phi'}) and a similar analysis to the above to obtain
\begin{equation}
\left| h_1(r)\right|\leq c\,r^{2-n}e^{I(r)}\max(\om(r),\sigma(r))\qquad\hbox{for}\ 0<r<1.
\label{eq:est-h1}
\end{equation}
Define $\zeta(r)$ by
\begin{equation}
\zeta(r)=\frac{h_1(r)}{h_0(r)}\quad\hbox{for}\ 0<r<1.
\label{def:zeta}
\end{equation}
Using  (\ref{def:h0}) -- (\ref{eq:est-h1}) we can estimate  $|\zeta(r)|, |r\zeta'(r)|\leq c\,\max(\om(r),\sigma(r))$. To estimate $\zeta''$, we write $h_0\zeta''=h_1''-h_0''\zeta-2 h_0'\zeta'$, where
\begin{equation}
h_0''(r)=(n-1)\,r^{-n}e^{I(r)}-r^{1-n}e^{I(r)}I'(r),
\label{eq:h0''}
\end{equation}
and
\begin{equation}\label{eq:h1''}
\begin{aligned}
h_1''(r)&=r^{-n}e^{I(r)}[(n-1)\tau(r)-rI'(r)\tau(r)-r\tau'(r)]\\
&+r^{-n}E_+(r)(n-1+R(r))\int_0^r \rho^{n-1}R(\rho)B[D^2v](\rho)d\rho-B[D^2v](r).
\end{aligned}
\end{equation}
 The terms $h_0''\zeta$ and $2h_0'\zeta'$ may be estimated pointwise as before, but $h_1''(r)$ involves the term $B[D^2v](r)$, which cannot be estimated pointwise. However, from (\ref{eq:Mp-B}) and $v\in X$  we conclude $M_p(r^2\zeta'',r)\leq c\,\max(\om(r),\sigma(r))$. Putting this together with the lower order derivatives, we obtain the desired estimate (\ref{eq:M2p(zeta)}). Summarizing so far, we have found a solution $Z$ of (\ref{eq:LZ=0}) in the desired form (\ref{eq:Z}).

Next we need to verify that any strong solution $u\in W_{\loc}^{2,p}(\overline{B_{1}}\backslash\{0\})$ of ${\mathcal L}u=0$ that satisfies the growth estimate (\ref{eq:u-hyp}) must be of the form 
(\ref{eq:u-asym}). To do this, we shall invoke well-known results for weighted Sobolev spaces. To begin with, let us introduce the weighted $L^p$-norm on $B_\circ=B_1\backslash\{0\}$:
\begin{equation}
\|u\|^p_{L^p_{\beta}(B_\circ)}=
\int_{0<|x|<1}|x|^{\beta p}\,|u(x)|^p\,dx.
\label{def:weightedLp(B*)}
\end{equation}
To relate this to the $M_p$-norm, notice that
\[
\int_{0<|x|<1}|x|^{\beta p}\,|u(x)|^p\,dx=\sum_{k=1}^\infty \int_{2^{-k}<|x|<2^{-k+1}}|x|^{\beta p}|u(x)|^p\,dx.
\]
Moreover, there exist constants $c_1$, $c_2$ (independent of $k$) such that
\[
c_1 2^{-k(\beta p+n)}\,M_p(u,2^{-k})^p\leq \int_{2^{-k}<|x|<2^{-k+1}}|x|^{\beta p}|u(x)|^p\,dx \leq c_2 2^{-k(\beta p+n)}\,M_p(u,2^{-k})^p.
\]
Consequently, $M_p(u,r)\leq c\,r^{\alpha}$ for $0<r<1$ implies $u\in L_\beta^p(B_\circ)$ if $\alpha+\beta>-n/p$, and conversely, $u\in L_\beta^p(B_\circ)$ implies $M_p(u,2^{-k})\leq c\,2^{k\beta+kn/p}$ which implies $M_p(u,r)\leq c_\alpha \,r^\alpha$ for $0<r<1$ if we choose $\alpha=-\beta-n/p$. We obtain analogous relationships between the $M_{p}$-norm and the $L^p_{\beta}$-norm of the terms $|x|^{|\alpha|}|\del^\alpha u(x)|$ for $|\alpha|\leq 2$.

Now let us introduce a weighted $L^p$-norm for functions on $\RR^{n}_\circ=\RR^n\backslash\{0\}$ with separate weights at the origin and infinity:
\begin{equation}\label{eq:weightedLp}
\begin{aligned}
\|u\|^p_{L^p_{\beta,\gamma}(\RR^{n}_\circ)}&=\|u\|^p_{L^p_{\beta}(B_\circ)}+
\|u\|^p_{L^p_{\gamma}(B^c)} \\
&=\int_{0<|x|<1}|u(x)|^p\,|x|^{\beta p}\,dx +\int_{|x|>1}|u(x)|^p\,|x|^{\gamma p}\,dx, 
\end{aligned}
\end{equation}
where $B^c=\RR^n\backslash \overline{B_1}$.
We then define the weighted Sobolev space $W^{2,p}_{\beta,\gamma}(\RR^{n}_\circ)$ to be those functions in $W^{2,p}_\loc(\RR^{n}_\circ)$ for which the norm
\begin{equation}
\|u\|_{W_{\beta,\gamma}^{2.p}}=
\sum_{|\alpha|\leq 2} \||x|^{|\alpha|}\partial^\alpha\, u\|_{L^p_{\beta,\gamma}(\RR^{n}_\circ)}
\label{eq:weightedSobolev}
\end{equation}
is finite.
Many authors have used similar weighted Sobolev spaces to study operators like the Laplacian
on $\RR^n$, $\RR^{n}_\circ$, and other noncompact manifolds with conical or cylindrical ends.
Using the analysis in \cite{MP},  \cite{Mc} or \cite{LM}, for example,
it is easily verified that the bounded operator
\begin{equation}
\Lap: W^{2,p}_{\beta,\gamma}(\RR^{n}_\circ)\to L^p_{\beta+2,\gamma+2}(\RR^{n}_\circ)
\label{eq:Lap-weighted}
\end{equation}
is Fredholm (finite nullity and finite deficiency) for all values of $\beta$ and $\gamma$ {\it except}
for the values $-2+\frac{n}{q}+k$ and $-\frac{n}{p}-k$ where $q=p/(p-1)$ and $k$ is any nonnegative integer. In fact, 
(\ref{eq:Lap-weighted}) is an isomorphism for $-n/p<\beta,\gamma<-2+n/q$ (recall that we are assuming $n\geq 3$, so such $\beta,\gamma$ exist). Since we are principally interested in the behavior of functions at the origin, we will fix $\gamma_0\in (-n/p,-2+n/q)$. Then
\[
\beta_-<\beta_+ \quad\Rightarrow \quad 
W_{\beta_-,\gamma_0}^{2,p}(\RR^{n}_\circ)\subset  W_{\beta_+,\gamma_0}^{2,p}(\RR^{n}_\circ).
\]
Moreover, for $\beta_+\in(-2+n/q,-1+n/q)$, we find that (\ref{eq:Lap-weighted}) is surjective with a one-dimensional nullspace spanned by $|x|^{2-n}\in W_{\beta_+,\gamma_0}^{2,p}(\RR^{n}_\circ)$; for $\beta_-\in(-n/p-2,-n/p-1)$, we find that (\ref{eq:Lap-weighted}) is injective with codimension equal to $n+1$.

Next we use perturbation theory (cf.\ \cite{K}, Ch.IV, Sec.5) to conclude not only that the operator
\begin{equation}
{\mathcal L}: W^{2,p}_{\beta,\gamma}(\RR^{n}_\circ)\to L^p_{\beta+2,\gamma+2}(\RR^{n}_\circ)
\label{eq:L-weighted}
\end{equation}
is Fredholm for exactly the same values of $\beta$ and $\gamma$ as for 
(\ref{eq:Lap-weighted}), but  the nullity and deficiency of 
(\ref{eq:Lap-weighted}) and (\ref{eq:L-weighted}) agree, provided $\de$ is sufficiently small.
So, in addition to the fixed $\gamma_0\in (-n/p,-2+n/q)$, let us now choose $\beta_0 \in (-n/p, -2+n/q)$, as well as $\beta_-\in(-1-n/p,-n/p)$ and $\beta_+\in(-2+n/q,-1+n/q)$.  Then
we assume that $\de$ is so small that
\begin{equation}
{\mathcal L}: W^{2,p}_{\beta_0,\gamma_0}(\RR^{n}_\circ)\to L^p_{\beta_0+2,\gamma_0+2}(\RR^{n}_\circ)
\quad\hbox{is an isomorphism},
\label{eq:L0-weighted}
\end{equation}
\begin{equation}
{\mathcal L}: W^{2,p}_{\beta_-,\gamma_0}(\RR^{n}_\circ)\to L^p_{\beta_- +2,\gamma_0+2}(\RR^{n}_\circ)
\quad\hbox{is injective with codim $= n+1$},
\label{eq:L--weighted}
\end{equation}
and
\begin{equation}
{\mathcal L}: W^{2,p}_{\beta_+,\gamma_0}(\RR^{n}_\circ)\to L^p_{\beta_+ +2,\gamma_0+2}(\RR^{n}_\circ)
\quad\hbox{is surjective with nullity $=1$}. 
\label{eq:L+-weighted}
\end{equation}
We claim that  $Z\in W^{2,p}_{\beta_+,\gamma_0}(\RR^{n}_\circ)$. In fact, this is quite simple to check given  the explicit formulas  (\ref{def:h0}) and (\ref{def:h1}), and the fact that $v\in X$, where $X$ has the norm (\ref{def:X}). Thus $Z$ is a basis vector for the one-dimensional nullspace of (\ref{eq:L+-weighted}). 

Now suppose $u\in W_{\loc}^{2,p}(\overline{B_{1}}\backslash\{0\})$ satisfies ${\mathcal L}u=0$ and the growth estimate (\ref{eq:u-hyp}) with $\e_0\in (0,1)$. Introduce a cut-off function $\chi\in C_0^\infty(B_1)$ equal to $1$ on $B_{1/2}$. If we now specify that $\beta_+\in(-1-\e_0+n/q,-1+n/q)$, then
  $\chi u\in W^{2,p}_{\beta_+,\gamma_0}(\RR^{n}_\circ)$. Let
$f={\mathcal L}(\chi u)$. Since $f=0$ for $|x|<1/2$ and for $|x|>1$, $f\in L^p_{\beta_0+2,\gamma_0+2}(\RR^{n}_\circ)$. But (\ref{eq:L0-weighted}) is an isomorphism, so we can find $v={\mathcal L}^{-1}f\in W^{2,p}_{\beta_0,\gamma_0}(\RR^{n}_\circ)$.
Now $\chi u-v\in W^{2,p}_{\beta_+,\gamma_0}(\RR^{n}_\circ)$ satisfies
${\mathcal L}(\chi u-v)=0$. Since the nullspace of (\ref{eq:L+-weighted}) is spanned by $Z$, there exists a constant $C$ such that $\chi u-v=CZ$.
But this means in particular that  $u=CZ+v$ for $0<|x|<1/2$. 

Now let us describe $v$ asymptotically. Let $\zeta_0,\zeta_1,\dots,\zeta_n$ denote a basis for the cokernel of (\ref{eq:L--weighted}), i.e.\  the $\zeta_i$ are linear functionals on
$L^p_{\beta_- +2,\gamma_0+2}(\RR^{n}_\circ)$ that are linearly independent and vanish on the image of (\ref{eq:L--weighted}). Now we want to choose $C_0,C_1,\dots,C_n$ so that 
${\mathcal L}(\chi(v-C_0-\sum_{j=1}^nC_jx_j))$ is in the range of (\ref{eq:L--weighted}), i.e.\ 
\begin{equation}
\zeta_i[{\mathcal L}(\chi(C_0+\sum_{j=1}^nC_jx_j)]=\zeta_i[{\mathcal L}(\chi v)]
\quad\hbox{for}\ i=0,\dots,n.
\label{eq:system1}
\end{equation}
To be able to solve the linear system (\ref{eq:system1}) for $C_0,\dots,C_n$, we need to verify that the finite-dimensional linear map
\begin{equation}
(C_0,\dots,C_n)\to \left(\zeta_0[{\mathcal L}(\chi(C_0+\sum_{j=1}^nC_jx_j)],\dots,
\zeta_n[{\mathcal L}(\chi(C_0+\sum_{j=1}^nC_jx_j)]\right)
\label{eq:linearmap}
\end{equation}
is nonsingular. Suppose $(C_0,\dots,C_n)$ is in the nullspace of (\ref{eq:linearmap}). Then ${\mathcal L}(\chi(C_0+\sum_{j=1}^nC_jx_j]$ is in the range of (\ref{eq:L--weighted}) and there exists  $\psi\in W^{2,p}_{\be_-,\gamma_0}(\RR^{n}_\circ)$
such that 
\[
{\mathcal L}(\psi)={\mathcal L}(\chi(C_0+\sum_{j=1}^nC_jx_j)),
\]
i.e.\ 
${\mathcal L}(\psi -\chi(C_0+\sum_{j=1}^nC_jx_j))=0$. But 
$\psi-\chi(C_0+\sum_{j=1}^nC_jx_j)\in W^{2,p}_{\beta_0,\gamma_0}(\RR^{n}_\circ)$ and (\ref{eq:L0-weighted}) is an isomorphism, so
$\psi=\chi(C_0+\sum_{j=1}^nC_jx_j)$. However, $\chi$ and $\chi x_j$ are {\it not} in $ W^{2,p}_{\be_-,\gamma_0}(\RR^{n}_\circ)$, so the only way that we can have $\psi=\chi(C_0+\sum_{j=1}^nC_jx_j)$ is to have $C_0=C_1=\cdots=C_n=0$. Consequently, the linear map (\ref{eq:linearmap}) is nonsingular.

Thus we can find $C_0,\dots,C_n$ and $w\in W^{2,p}_{\be_-,\gamma_0}(\RR^{n}_\circ)$ so that 
${\mathcal L}(w)={\mathcal L}(\chi(v-C_0-C_1x_1\cdots C_nx_n))$. But
(\ref{eq:L0-weighted}) is an isomorphism, so we obtain
$v=C_0+C_1x_1+\cdots+C_nx_n+w$ for $0<|x|<1/2$.
This yields (\ref{eq:u-asym}) and we only need to verify $M_{2,p}(w,r)\leq c\,r^{2-\e_1}$.
But recall that $w\in W^{2,p}_{\be_-,\gamma_0}(\RR^{n}_\circ)$
implies $M_{2,p}(w,r)\leq c\,r^{-\beta_- -n/p}$ for $0<r<1/2$, and we can then let
$\beta_-=-2-n/p+\e_1$ for any $\e_1\in (0,1)$ to obtain the desired estimate.
\hfill $\Box$

\medskip
Now let us formulate the result for a general point $y\in\RR^n$ where we do not assume $a_{ij}(y)=\de_{ij}$. With $y$ fixed and the same conditions (\ref{eq:om1}) and (\ref{eq:om2}) on $\om$,  we now assume
\begin{equation}
\sup_{ |x-y|=r}\|{\bf A}_x-{\bf A}_y\|\leq
\om(r)\quad\hbox{for}\ 0<r<1,
\label{eq:aij-asym}
\end{equation}
and we want to construct a solution of 
\begin{equation}
{\mathcal L}(x,\del_x)Z_y(x)=0 \quad\hbox{for}\  x\in B_\e(y)\backslash\{y\},
\label{eq:LZ_y=0_2}
\end{equation}
 for $\e$ sufficiently small.
(Unlike Theorem 1, in the Corollary below we need to assume that the coefficients are real-valued so that we can choose coordinates in which $a_{ij}(0)=\de_{ij}$.)

\begin{cor}
For $n\geq 3$ and $p\in (1,\infty)$, fix $y\in\RR^n$ and suppose that the constant coefficient operator ${\mathcal L}(y,\del_x)$ is elliptic and the coefficients $a_{ij}(x)$ are bounded, measurable functions satisfying (\ref{eq:aij-asym}).  For  $\e>0$ sufficiently small, there exists a solution  of (\ref{eq:LZ_y=0_2}) in the form 
\begin{equation}
Z_y(x)=h_y(|{\bf A}_y^{-1/2}(x-y)|)+v(x)
\label{eq:Zy}
\end{equation}
 where $h_y$ is defined by 
\begin{equation}
h_y(r)=\int_{r}^{\e} s^{1-n}\, e^{I_y(s)}\,ds \,(1+\zeta_y(r)),
\label{eq:h_y}
\end{equation}
with $I_y(r)$ given by (\ref{def:I_y}) and 
\begin{equation}
M_{2,p}(\zeta_y,r;y)\leq c\,\max\left(\om(r),\sigma(r)\right),
\label{eq:M2p(zeta_y)}
\end{equation}
and $v$ in (\ref{eq:Zy}) satisfies
\begin{equation}
M_{2,p}(v,r;y)\leq c\,r^{2-n}\,e^{I_y(r)}\,\om(r).
\label{eq:M2p(v,y)} 
\end{equation}
Moreover, for any  $u\in W_{\loc}^{2,p}(\overline{B_{\e}(y)}\backslash\{y\})$ that is a strong
solution of ${\mathcal L}(x,\del_x)u=0$ in $\overline{B_{\e}(y)}\backslash\{y\}$ subject to the growth condition 
\[
M_{2,p}(u,r;y)\leq c\, r^{1-n+\e_0}
\] where $\e_0>0$,
 there exist constants $C, C_0, C_1,\dots,C_n$ (depending on $u$) such that
\begin{equation}
u(x)=C\,Z_y(x)+C_0+\sum_{j=1}^n C_j(x_j-y_j)+w(x) \quad\hbox{for}\ 0<|x-y|<\e,
\label{eq:u-asym2}
\end{equation}
where $w$ satisfies
\begin{equation}
M_{2,p}(w,r;y)\leq c\,r^{2-\e_1} \quad\hbox{for any}\ \e_1>0.
\label{eq:M2p(w)2}
\end{equation}
\end{cor}

If we use integration by parts, we can  write the solution in Corollary 1 as
\begin{equation}
Z_y(x)=\frac{\langle{\bf A}_y^{-1}(x-y),(x-y)\rangle^{\frac{2-n}{2}}}{(n-2)}\,
e^{I_y\left(\sqrt{\langle {\bf A}_y^{-1}(x-y),(x-y)\rangle}\right)}\,(1+\xi_y(x)),
\label{eq:Zy2}
\end{equation}
where 
$M_{1,\infty}(\xi_y,r;y)\leq c\,\max(\om(r),\sigma(r))$ for $0<r<\e$.
Notice that, if $I_y(r)$ has a finite limit as $r\to 0$, then the leading term in 
 (\ref{eq:Zy2}) is just a constant times $\tilde F_y(x-y)$, the fundamental solution for ${\mathcal L}(y,\del_x)$ at $y$ (cf.\  (\ref{def:F_y})).

\smallskip
\noindent{\bf Proof of Corollary 2.} For now we continue to assume $y=0$, i.e.
\begin{equation}
\sup_{ |x|=r}\|{\bf A}_x-{\bf A}_0\|\leq
\om(r),
\label{eq:aij-asym0}
\end{equation}
and let us assume ${\bf A}_0$ is positive definite.
 Let ${\bf B}=(b_{ij})={\bf A}_0^{-1/2}$ so that
${\bf B} {\bf A}_0 {\bf B} ={\bf I}.$
Introduce new independent variables $\tilde{x}={\bf B} x$ and the matrix
$\tilde{\bf A}_{\tilde x}=(\tilde a_{ij}(\tilde x))= {\bf B}{\bf A}_x{\bf B}$,
which satisfies $\tilde a_{ij}(0)=\de_{ij}$. Using $\del/\del x_i=b_{ik}\del/\del \tilde x_k=b_{ki}\del/\del \tilde x_k$, we can then write
\begin{equation}
{\mathcal L}(x,\del_x)=a_{ij}(x)\frac{\del^2}{\del x_i\del x_j}
=\tilde a_{k\ell}(\tilde x)\frac{\del^2}{\del \tilde x_k\del \tilde x_\ell}=\tilde{\mathcal L}(\tilde x,\del_{\tilde x})
\label{eq:L=tilde-L}
\end{equation}
and apply Theorem 1 to $\tilde{\mathcal L}$ in the coordinates $\tilde x$. We conclude that for $\e>0$ sufficiently small, there is a solution $\tilde Z$ of
$\tilde{\mathcal L}(\tilde x,\del_{\tilde x})\tilde Z( \tilde x)=0$ for $0<|\tilde x|<\e$,
of the form $\tilde Z(\tilde x)=\tilde h(|\tilde x|)+\tilde v(\tilde x)$ where $M_{2,p}(\tilde v,r)\leq c\,r^{2-n}e^{I(r)}\om(r)$ and 
$\tilde h(r)$ is of the form (\ref{est:h}) with
\begin{equation}
I(r)=\frac{1}{|S^{n-1}|}\int_{r<|\tilde z|<\e}
\left({\rm tr}\,\tilde{\bf A}_{ \tilde z}-n\langle\tilde{\bf A}_{ \tilde z}  \tilde z,  \tilde z\rangle | \tilde z|^{-2}\right)\,\frac{d \tilde z}{|\tilde z|^n}.
\label{eq:I(r)-tilde}
\end{equation}
Expressed in terms of the original variables $x$, we obtain
$Z(x)=\tilde Z(\tilde x)=\tilde h(|{\bf B}x|)+\tilde v({\bf B}x)$ that satisfies
${\mathcal L}( x,\del_{ x})Z( x)=0$ for $0<|{\bf B} x|<\e$; but choosing $\e_0$ sufficiently small, we conclude that ${\mathcal L}( x,\del_{ x})Z( x)=0$ for $0<| x|<\e_0$.

Finally, if $y$ is a general point in $\RR^n$, then let $\tilde x={\bf B}(x-y)$ with ${\bf B}={\bf A}_y^{-1/2}$
and let $\tilde{\bf A}_{\tilde x}=(\tilde a_{ij}(\tilde x))= {\bf B}{\bf A}_x{\bf B}$; since $\tilde x=0$ corresponds to $x=y$, we have $\tilde a_{ij}(0)=\de_{ij}$ and we can apply Theorem 1 to 
$\tilde{\mathcal L}(\tilde x,\del_{\tilde x})={\mathcal L}( x,\del_{ x})$.
We  obtain the solution $Z(x)=\tilde h(|{\bf B}(x-y)|)+\tilde v({\bf B}(x-y))$
where $\tilde h(r)$ involves $I(r)$ as in (\ref{eq:I(r)-tilde}). To transform (\ref{eq:I(r)-tilde}) to the original variables, replace $\tilde{\bf A}_{\tilde z}$ by ${\bf A}_z$ and every other occurrence of $\tilde z$ by $z-y$; we find that (\ref{eq:I(r)-tilde}) is of the desired form  (\ref{def:I_y}),
so we may let $h_y(r)=\tilde h(r)$. Moreover, since $\tilde v$ satisfies $M_{2,p}(\tilde v,r)\leq c\,r^{2-n}e^{I(r)}\om(r)$, it is clear that $v(x)=\tilde v({\bf B}(x-y))$ satisfies $M_{2,p}(v,r;y)\leq c\,r^{2-n}e^{I(r)}\om(r)$, completing the proof.
\hfill $\Box$

\section{Finding the Constant $C_y$ so that ${\mathcal L}Z_y(x)=C_y\delta(x-y)$ in $B_\e(y)$}

Let us now begin to discuss the role that $Z_y(x)$ plays in finding the fundamental solution. As before, we first consider $y=0$ with $a_{ij}(0)=\de_{ij}$ and then use a change of variables to consider a general fixed $y\in\RR^n$; as in the proof of Theorem 1, we shall assume $\e=1$.
We first want to see whether the function $Z(x)$ found in Theorem 1 satisfies
\begin{equation}
-{\mathcal L}(x,\del_x)Z(x)=C_0\delta(x)\quad\hbox{ for} \ x\in B_1(0)
\label{eq:-LZ=delta0}
\end{equation}
for some constant $C_0$. 

It is not immediately clear how the left-hand side of (\ref{eq:-LZ=delta0}) is defined. Recall from the proof of Theorem 1 the decomposition $Z(x)=h(|x|)+v(x)$, where $h$ is given by (\ref{eq:h})  and $M_{2,p}(v,r)\leq c\,r^{2-n}\,e^{I(r)}\,\om(r)$. We can easily calculate $\del_i\del_j Z$ and show that for any $\mu>0$ there is a constant $C_\mu$ so that 
\begin{equation}
M_p(\del_i\del_j Z,r)\leq  C_\mu \, r^{-n-\mu}.
\label{est:Mp(D^2Z)}
\end{equation}
In fact, since the $a_{ij}$ are bounded functions, we conclude that for each $i,j$ the function
\begin{equation}
F_{ij}(x) =a_{ij}(x)\del_i\del_j Z(x) \quad\hbox{for}\ x\not=0
\label{def:Fij}
\end{equation}
satisfies
\begin{equation}
M_p(F_{ij},r)\leq  C_\mu \, r^{-n-\mu} \quad\hbox{for}\ 0<r<1.
\label{est:Mp(aD^2Z)}
\end{equation}
But this estimate implies that $F_{ij}$ can be regularized at $x=0$ to give a distribution ${\mathcal F}_{ij}$, in particular as a continuous linear functional on the space of $\lambda$-H\"older continuous functions of compact support in $U=B_1(0)$:
\begin{equation}
\langle{\mathcal F}_{ij},\phi\rangle=\int_{|x|<1} F_{ij}(x)\left[ \phi(x)-\phi(0)\chi(|x|)\right]\,dx
\label{def:functionalFij}
\end{equation}
where $\chi(r)$ is a smooth cut-off function which is identically 1 near $r=0$ but vanishes for $r>1/2$.
(Since $|\phi(x)-\phi(0)\chi(|x|)|\leq C|x|^\lambda$, by choosing $0<\mu<\lambda$ we see that the integral in (\ref{def:functionalFij}) converges.)
Now let us take the sum over all $i,j$:
\begin{equation}
{\mathcal F}_0=\sum_{i,j}\,{\mathcal F}_{ij}.
\label{def:functionalF0}
\end{equation}
Then ${\mathcal F}_0$ is a regularization of ${\mathcal L}(x,\del_x)Z(x)=0$, and so it vanishes:
\[
\langle {\mathcal F}_0,\phi\rangle = \sum_{i,j}\langle {\mathcal F}_{ij},\phi\rangle
=\int \sum_{i,j}F_{ij}(x)\left[ \phi(x)-\phi(0)\chi(|x|)\right]\,dx =0.
\]
 Of course, regularization effects the distribution only at $x=0$, so {\it if} we can interpret the expression
\begin{equation}
{\mathcal F}={\mathcal L}(x,\del_x)Z(x),
\label{def:functionalF}
\end{equation}
as a distribution, then ${\mathcal F}$  is supported only at $x=0$. As such, it is a linear combination of the delta distribution and its derivatives. But since ${\mathcal F}$ is a continuous linear functional on 
$\lambda$-H\"older functions with $\lambda\in (0,1)$, it must only involve the delta distribution itself, i.e.\  (\ref{eq:-LZ=delta0}) must hold for some constant $C_0$. But we still have two problems: 
(i) how is ${\mathcal F}$ itself defined as a distribution? and (ii) how do we calculate the constant $C_0$?

The difficulty in defining ${\mathcal F}$ as a distribution in $U$ is a consequence of the lack of regularity of the $a_{ij}$. In particular, there is no difficulty in defining the distributional derivatives of $Z$:
\begin{equation}
\langle \del_i\del_j Z, \phi\rangle =-\int_{U} \del_j Z(x)\del_i\phi(x)\,dx
\label{eq:didjZ}
\end{equation}
for $\phi\in C_0^1(U)$, since the integral on the right converges. To handle the $a_{ij}$, we can write\begin{equation}
\langle {\mathcal L}Z,\phi\rangle=\int_{U} \left[(a_{ij}(x)-\delta_{ij})\,\del_i\del_j Z(x)\,\phi(x)
- \del_iZ(x)\del_i\phi(x)\right]\,dx.
\label{eq:<LZ,phi>}
\end{equation}
Of course, the integral in (\ref{eq:<LZ,phi>}) is actually an improper integral due to the singularity 
of $\del_i\del_j Z$ at $x=0$; but 
provided this integral converges, we conclude that $-{\mathcal L}Z=C_0\de$ in $U$, and
 we can calculate $C_0$  by
\begin{equation}
C_0=\lim_{\e\to 0}\,\int_{|x|<1} \left[(-a_{ij}(x)+\delta_{ij})\,\del_i\del_j Z(x)\,\phi_\e(|x|)
+ \del_iZ(x)\del_i\phi_\e(|x|)\right]dx,
\end{equation}
where $\phi_\e(|x|)=\chi(|x|/\e)$ with the cut-off function $\chi$  introduced above;
for these purposes, we are able to assume $\phi(x)=\phi(|x|)$ since we can write
$\phi(x)=\phi_0(|x|)+\phi_1(x)$ with $|\phi_1(x)|+|x|\,|\nabla\phi_1(x)|\leq c\,|x|$ for $|x|<1$, which shows that
$\langle {\mathcal L}Z,\phi_1\rangle$ is well-defined and $\phi_1$ contributes nothing to $C_0$.
We shall now prove the following.

\begin{theorem} Suppose the conditions of Theorem 1 hold and $Z$ is the function found there.
If $I(0)=\lim_{r\to 0}I(r)$ exists and is finite then the integral in (\ref{eq:<LZ,phi>}) converges and we can calculate $C_0=|S^{n-1}|\,e^{I(0)}$.
If $I(r)\to -\infty$ as $r\to 0$, then the integral in (\ref{eq:<LZ,phi>}) converges and we can calculate $C_0=0$.
\end{theorem}

\noindent
{\bf Proof of Theorem 2.}
Recall from the proof of Theorem 1 the decomposition $Z(x)=h(|x|)+v(x)$, where $h$ is given by (\ref{eq:h})  and $v$ satisfies (\ref{eq:M2p(v)}). Since we always assume that $I(r)$ is bounded above, we obtain from (\ref{eq:M2p(v)})
\begin{equation}
M_{2,p}(v,r)\leq c\,r^{2-n}\,\om(r) \quad\hbox{for}\ 0<r<1.
\label{eq:om-2}
\end{equation}
We shall separately consider the roles of $v$ and $h$. In the estimates below, $\int_{|x|<\e}$ should actually be considered as an improper integral $\lim_{\eta\to 0}\int_{\eta<|x|<\e}$, but we avoid such cumbersome notation.

For $v$, we use (\ref{est:int-Mp}) and (\ref{eq:om-2}) 
to conclude 
\[
\left|\int_{|x|<\e}\del_i v\,\del_i\phi_\e\,dx\right|
\leq c\,\e^{-1}\int_0^{\e}\om(r)\,dr \leq c\,\om(\e)
\to 0.
\]
 and 
\[
\left|\int_{|x|<\e}(a_{ij}-\de_{ij})\,\del_i\del_jv\,\phi_\e\,dx\right|
\leq c\,\int_0^{\e}\om^2(r)\,\frac{dr}{r} =c\,\sigma(\e) \to 0.
\]
This shows that $v$ makes no contribution to the value of $C_0$.

 To determine the effect of $h$, let us write $h(r)=h_2(r)+h_3(r)+c$ where
 \begin{equation}
 h_2(r)=c_1\int_r^1s^{1-n}E_+(s)\,ds
 \end{equation}
 with $c_1$ chosen as in the proof of Theorem 1 so that $c_1\,E_+(r)=e^{I(r)}+o(1)$ as $r\to 0$, and 
 \begin{equation}
 h_3(r)=\int_r^1 s^{1-n}E_+(s)\int_0^s\rho^{n-1}E_-(\rho)B[D^2 v](\rho)\,d\rho.
 \end{equation}
 (Notice that $h_2$ and $h_3$ differ slightly from (\ref{def:h0}) and (\ref{def:h1}).)
 Let us consider $h_3$ first:
\[
\del_i h_3=-x_ir^{-n}E_+(r)\int_0^r \rho^{n-1}E_-(\rho)B[D^2v](\rho)\,d\rho
\]
and
\[
\del_i\del_j h_3=
-r^{-n}E_+(r)\left(\de_{ij}-n\frac{x_ix_j}{r^2}+\frac{x_ix_j}{r^2}R(r)\right)\int_0^r\rho^{n-1}E_-(\rho)B[D^2v](\rho)\,d\rho
\]
\[
-\frac{x_ix_j}{r^2}B[D^2v](r).
\]
The calculation of $\del_ih_3$ combined with 
 (\ref{eq:int-phi'}) and  the assumption that $I(r)$ is bounded  shows
\begin{equation}
\left|\int_{|x|<\e}\del_i h_3\,\del_i\phi_\e\,dx\right|
\leq 
c\,\e^{-1}\int_0^{\e}\sigma(r) \,dr\leq c\,\sigma(\e)\to 0.
\end{equation}
On the other hand, the ``worst'' term in $(a_{ij}-\de_{ij})\del_i\del_j h_3$ is
\[
W(x)=-r^{-n}E_+(r)\left(a_{ii}-n\frac{a_{ij}x_ix_j}{r^2}-\left(\frac{a_{ij}x_ix_j}{r^2}-1\right)R\right)\int_0^r\rho^{n-1}E_-(\rho)B[D^2v](\rho)\,d\rho.
\]
Using the fact that $\phi_\e$ only depends on $r$ and observing that
\[
\int_{S^{n-1}}\left(a_{ii}-n\frac{a_{ij}x_ix_j}{r^2}-\left(\frac{a_{ij}x_ix_j}{r^2}-1\right)R\right)\,d\theta=|S^{n-1}|R(r),
\]
 we can calculate
\[
\langle W,\phi_\e\rangle = -|S^{n-1}|\int_0^\e E_+(r)\frac{R(r)}{r}\,\phi_\e(r)\int_0^r\rho^{n-1}E_-(\rho)B[D^2v](\rho)\,d\rho\,dr,
\]
\[
=|S^{n-1}|\int_0^\e (E_+)'(r)\,\phi_\e(r)\int_0^r\rho^{n-1}E_-(\rho)B[D^2v](\rho)\,d\rho\,dr.
\]
Now we can integrate by parts to obtain
\[
\langle W,\phi_\e\rangle =
 |S^{n-1}|\left(\int_0^\e E_+(r)\,\phi'_\e(r)\int_0^r\rho^{n-1}E_-(\rho)B[D^2v](\rho)\,d\rho\,dr\right.
\]
\[
\left.+\int_0^\e \,\phi_\e(r)\,r^{n-1}B[D^2v](r)\,dr\right).
\]
But, using $|\phi'_\e|\leq c\,\e^{-1}$ and $|E_+|\leq c$, we find
\[
\left|\int_0^\e E_+(r)\,\phi'_\e(r)\int_0^r\rho^{n-1}E_-(\rho)B[D^2v](\rho)\,d\rho\,dr\right|
\leq c\,\e^{-1}\,\int_0^\e \sigma(r)\,dr\leq c\,\sigma(\e) \to 0,
\]
and
\[
\left|
\int_0^\e \phi_\e(r)\,r^{n-1}B[D^2v](r)\,dr
\right|
\leq
c\,\sigma(\e)\to 0.
\]
We conclude that $h_3$ makes no contribution to $C_0$.

Now consider $h_2$. We compute
$\del_i h_2=-c_1 x_i r^{-n}E_+(r)$ 
and 
\[
\del_i\del_j h_2=-c_1\,r^{-n}\,E_+(r)\left(\de_{ij}-n\frac{x_ix_j}{r^2}-\frac{x_ix_j}{r^2}R(r)\right),
\]
so
\[
(a_{ij}-\de_{ij})\del_i\del_j h_2 = -c_1\, r^{-n}\, E_+(r) \left(a_{ii}-n\frac{a_{ij}x_ix_j}{r^2}-\left(\frac{a_{ij}x_ix_j}{r^2}-1\right)R\right).
\]
Notice that
\[
\int_{|x|<\e} r^{-n}\, E_+(r) \left(a_{ii}-n\frac{a_{ij}x_ix_j}{r^2}-\left(\frac{a_{ij}x_ix_j}{r^2}-1\right)R\right)\chi(\frac{r}{\e})\,dx
\]
\[
=|S^{n-1}|\int_0^\e E_+(r)\frac{R(r)}{r}\chi(\frac{r}{\e})\,dr,
\]
so we now may compute
\[
-\langle {\mathcal L}h_2,\phi_\e\rangle =
c_1|S^{n-1}|\int_0^\e \left( E_+(r)\frac{R(r)}{r}\chi(\frac{r}{\e})-E_+(r)\left[\chi(\frac{r}{\e})\right]'\right)\,dr
\]
\[
=-c_1|S^{n-1}|\int_0^\e \left[E_+(r)\chi(\frac{r}{\e})\right]'\,dr.
\]
Now, if $I(0)$ exists and is finite, then $c_1E_+(0)=e^{I(0)}$, and we conclude
$-\langle {\mathcal L}h_2,\phi_\e\rangle =|S^{n-1}|e^{I(0)}$. On the other hand, 
if $I(r)\to -\infty$ as $r\to 0$, then $E_+(r)\to 0$, so $C_0=0$.
This completes the proof of Theorem 2.
 \hfill $\Box$

  \medskip
 Now let us consider a general fixed $y\in \RR^n$ and try to find the constant $C_y$ so that
 \begin{equation}
-{\mathcal L}(x,\del_x)Z_y(x)=C_y\delta(x)\quad\hbox{ for} \ x\in B_1(y).
\label{eq:-LZ=delta}
\end{equation}
We replace (\ref{eq:<LZ,phi>}) by
\begin{equation}
\langle {\mathcal L}Z_y,\phi\rangle=\int_{U} \left[(a_{ij}(x)-a_{ij}(y))\,\del_i\del_j Z_y(x)\,\phi(x)
- a_{ij}(y)\del_iZ_y(x)\del_j\phi(x)\right]\,dx.
\label{eq:<LZy,phi>}
\end{equation}

 \begin{cor}
Suppose the conditions of Corollary 1 hold and $Z_y$ is the function found there.
 If  $I_y(0)=\lim_{r\to 0}I_y(r)$ exists and is finite then the integral in (\ref{eq:<LZy,phi>}) converges and $Z_y$ satisfies (\ref{eq:-LZ=delta}) with $C_y$ given by 
 \begin{equation}
 C_y=|S^{n-1}|\sqrt{{\rm det}{\bf A}_y}\,e^{I_y(0)}.
 \label{def:Cy}
 \end{equation}
If $I_y(r)\to -\infty$ as $r\to 0$ then the integral in (\ref{eq:<LZy,phi>}) converges and $Z_y$ satisfies (\ref{eq:-LZ=delta}) with $C_y=0$.
\end{cor}

\medskip\noindent
{\bf Proof of Corollary 3.}
We need only show $\langle -{\mathcal L}Z_y,\phi\rangle =|S^{n-1}|(\hbox{det}{\bf A}_y)^{1/2}e^{I_y(0)}\phi(y)$ for some $\phi\in C_0^\infty(B_{\e_y}(y))$. Let us recall the change of coordinates used in the proof of Theorem 1, namely $\tilde x={\bf B}(x-y)$ where ${\bf B}={\bf A}_y^{-1/2}$, and let
 $\tilde\phi(\tilde x)=\phi(x)$; then
 \[
-\int {\mathcal L}(x,\del_x)Z_y(x)\phi(x)\,dx
=-(\hbox{det}{\bf A}_y)^{1/2}\int \tilde{\mathcal L}(\tilde x,\del_{\tilde x})\tilde Z_0(\tilde x)\tilde\phi(\tilde x)\,d\tilde x.
\]
But  Theorem 2 implies $-\langle \tilde{\mathcal L} \tilde Z_0, \tilde\phi\rangle=
|S^{n-1}|e^{I_y(0)}\tilde \phi(0)$. Since $\tilde\phi(0)=\phi(y)$, we obtain the desired result.  \hfill $\Box$

\section{Constructing the Fundamental Solution}

\medskip
Now we are in a position to construct the fundamental solution in a bounded open set $U\subset \RR^n$ when $a_{ij}\in C^\om(U)$, assuming that $\om$ satisfies the square-Dini condition (\ref{eq:om1}) and for every $y\in U$ we know that $I_y(0)$ exists and is finite. Given the results in the preceding section, it is natural to seek the fundamental solution in the form $F(x,y)=Z_y(x)/C_y+v(x,y)$ where ${\mathcal L}(x,\del_x)v(x,y)=0$. In fact, it will be convenient to construct $F(x,y)$ as 
 the  Green's function $G(x,y)$, in the sense of (\ref{def:Green's}), for a smooth, bounded domain $V$ that contains $U$.
 
 But first let us observe that the additional assumption that $I_y(0)$ exists and is finite allows us to improve the asymptotic description of $Z_y(x)$. In fact, let us fix $y=0$ with $a_{ij}(0)=\de_{ij}$, and assume
 \begin{equation}
 |I(r)-I(0)|\leq \theta(r)
 \label{eq:theta-0}
 \end{equation}
 where $\theta(r)$
is a positive, nondecreasing function for $0< r< 1$ such that $\theta(0)=0$; as with $\om$ we additionally assume that for some $\nu\in (0,1)$ we have
\begin{equation}
\theta(r)\,r^{-1+\nu}\quad\hbox{is nonincreasing for}\ 0<r<1.
\label{eq:theta2}
\end{equation}
 (If $\om(r)$ satisfies the Dini condition, then we can take $\theta(r)=\omega(r)$.)

\begin{lemma}
Under the conditions of Theorem 1, let us additionally assume that (\ref{eq:theta-0}) and (\ref{eq:theta2}) hold. Then the solution $Z$ found in Theorem 1 satisfies
\begin{equation}
Z(x)=\frac{|x|^{2-n}e^{I(0)}}{n-2}\left(1+\xi(x)\right),
\label{eq:Z-asym3}
\end{equation}
where for any $p\in (1,\infty)$ we can estimate $\xi$ by
\begin{equation}
M_{2,p}(\xi,r)\leq c\,\max(\om(r),\sigma(r),\theta(r)) \quad\hbox{as}\ r\to 0,
\label{eq:M2p(xi)}
\end{equation}
where $c$ depends only on $\om$, $\theta$, $n$, and ${I(0)}$.
\end{lemma}

\noindent
{\bf Proof of Lemma 1.} In the decomposition $Z(x)=h(|x|)+v(x)$ as in (\ref{eq:Z}), it is elementary to use (\ref{eq:M2p(v)}) with (\ref{eq:theta-0}) to show that (\ref{eq:M2p(xi)}) applies to $\xi=c\,|x|^{n-2}\,v$, so we focus on $h$. Recall the decomposition $h=h_0+h_1$ using (\ref{def:h0}) and (\ref{def:h1}). Let us recall $h_0''$ from (\ref{eq:h0''}) and use (\ref{eq:theta-0}) to estimate $e^{I(r)}\leq (1+2\theta(r))e^{I(0)}$ for $r$ sufficiently small; also recalling  $|I'(r)|\leq 2(n-1)\om(r)/r$, we obtain
\[
|h_0''(r)-(n-1)r^{-n}e^{I(0)}|
\leq c\, r^{-n}\,\max(\omega(r),\theta(r));
\]
here (and henceforth) $c$ depends only on $\om$, $\theta$, $n$, and $I(0)$.
 Using (\ref{eq:om2}) and (\ref{eq:theta2}), we can integrate this twice to obtain
\[
\left| h_0(r)-\frac{r^{2-n}e^{I(0)}}{n-2}\right|\leq c\,r^{2-n}\,\max(\omega(r),\theta(r)).
\]
In fact, this argument has shown that $\xi_0(r)=(n-2)\,r^{n-2}\,e^{-I(0)}\,h_0(r)-1$ satisfies the pointwise estimate
\[
M_{2,\infty}(\xi_0,r)\leq c\,\max(\om(r),\theta(r)).
\]
Similarly, let us recall $h_1''$ from (\ref{eq:h1''}); we can estimate most of the terms pointwise, but the term $B[D^2v](r)$ can only be estimated in $M_p$, so $\xi_1=(n-2)\,r^{n-2}\,e^{-I(0)}\,h_1(r)$ satisfies
\[
M_{2,p}(\xi_1,r)\leq c\,\max(\om(r),\sigma(r)).
\]
These may be combined to yield (\ref{eq:M2p(xi)}), so the Lemma is proved. \hfill $\Box$

\medskip
Now we are ready to construct our fundamental solution in a bounded open set.

\begin{theorem}
Suppose (\ref{eq:L}) is uniformly elliptic in a bounded open set $U\subset{\RR}^n$ for $n\geq 3$, where the coefficients 
$a_{ij}$ are continuous functions with uniform modulus of continuity $\om$ satisfying (\ref{eq:om1}).  If $I_y(0)$ exists at every $y\in U$ and 
\begin{equation}
|I_y(r)-I_y(0)|\leq \theta(r) \quad \hbox{for all} \ y\in U,
\label{eq:theta_y}
\end{equation}
where $\theta$ is a positive, nondecreasing function for $0< r< 1$ with $\theta(0)=0$ and 
(\ref{eq:theta2}),
then there is a function $F(x,y)$ satisfying (\ref{def:FS})  for $x,y\in U$; moreover, $F(x,y)$ admits the asymptotic description
\begin{equation}
F(x,y)= \frac{\langle{\bf A}_y^{-1}(x-y),(x-y)\rangle^{\frac{2-n}{2}}}{(n-2)|S^{n-1}|\,\sqrt{\det{\bf A}_y}}\,(1+H(x,y)),
\label{eq:G-asym0}
\end{equation}
where for any $p\in (1,\infty)$ and any compact set $K\subset U$ we have 
\begin{equation}
M_{2,p}(H(\cdot,y),r;y)\leq c\,\max(\om(r),\sigma(r),\theta(r)) \quad\hbox{as}\ r\to 0,
\label{eq:M2p(H)} 
\end{equation}
with constant $c$ independent of $y\in K$.
\label{th:3}
\end{theorem}

\noindent
{\bf Proof of Theorem 3.} For each $y\in U$, denote the $\e$ in Corollary 1 
by $\e_y$. Note that the size of $\e_y$  depends on the behavior of the coefficients $a_{ij}$ near $y$ through their ellipticity (i.e.\ the norm of ${\bf A}_y^{-1/2}$) and their continuity  (i.e.\ $\omega$). Since the ellipticity and modulus of continuity are uniform on $U$, we can find $\e>0$ that is independent of  $y\in U$. In fact, if we choose a smooth, bounded domain $V$ with $V\supset\overline U$   and $\hbox{dist}(U,\del V)>\e$, then we can extend the coefficients $a_{ij}$ to $V$ in such as way as to maintain the uniform ellipticity as well as the modulus of continuity $\om(r)$ (cf.\ \cite{Mi}).
Thus for all $y\in \overline U$ we can construct $Z_y(x)$ in $B_\e(y)\backslash\{y\}$  . In fact, repeating this argument with a neighborhood of $V$, we may assume that $Z_y(x)$ is defined for all $y\in \overline V$.

For each $y\in \overline V$ let us use Corollary 2 to calculate $C_y>0$, and conclude that
$- {\mathcal L}(x,\del_x)Z_y(x)/C_y=\de(x-y)$ for all $x,y\in \overline V$ with $|x-y|<\e$.  We shall construct $G(x,y)$ as the Green's function for ${\mathcal L}$ in $V$. For fixed $y\in V$, 
let us introduce a smooth cut-off function $\eta_y(r)$ satisfying $\eta_y(r)=1$ for sufficiently small $r>0$ but $\eta_y(|\cdot-y|)$ has compact support in $V$. Then let us define
\begin{equation}
G(x,y)=\eta_y(|x-y|)Z_y(x)/C_y + v(x,y),
\label{def:G(x,y)}
\end{equation}
where $v(x,y)$ is to be determined. But if we apply $-{\mathcal L}(x,\del_x)$ to $G(x,y)$ we obtain
\[
-{\mathcal L}(x,\del_x)G(x,y)=\delta(x-y)+\psi(x,y)-{\mathcal L}(x,\del_x)v(x,y),
\]
where $\psi(\cdot,y)\in L^p(V)$. So, for fixed $y\in V$, consider the Dirichlet problem for $v$:
\begin{equation}
\begin{aligned}
&{\mathcal L}(x,\del_x)v(x,y)=\psi(x,y) \quad\hbox{for}\  x\in V, \\ 
&v(x,y)=0 \quad\hbox{for}\  x\in \del V.
\end{aligned}
\label{def:v}
\end{equation}
It is well-known (cf.\ Theorem 9.15 in \cite{GT}) that (\ref{def:v}) has a unique solution $v(\cdot,y)\in W^{2,p}(V)\cap  W_0^{1,p}(V)$, so using this $v(x,y)$ in (\ref{def:G(x,y)})
not only ensures that (\ref{def:Green's}) holds for $x,y\in V$ (and hence for $x,y\in U$), but the following: for any $f\in C(\overline V)$,
$u(x)=-\int_V G(x,y)f(y)\,dy$ satisfies ${\mathcal L}u=f$ in $V$ and $u=0$ on $\del V$. Now if we pick
$\phi\in C^2(\overline V)$ with $\phi=0$ on $\del V$ and let $f={\mathcal L}\phi$, then uniqueness of the solution shows $u=\phi$ in $V$, i.e.\ (\ref{def:Green's}) holds. In other words, the $G(x,y)$ that we have constructed is just the Green's function for ${\mathcal L}$ in $V$; in particular, (\ref{def:FS}) holds.

Now from (\ref{def:G(x,y)}) and Corollaries 1 and 2 we see that
\begin{equation}
G(x,y)=\frac{\langle{\bf A}_y^{-1}(x-y),(x-y)\rangle^{\frac{2-n}{2}}}{(n-2)|S^{n-1}|\,\sqrt{\det{\bf A}_y}}\,e^{J_y(\sqrt{\langle {\bf A}_y^{-1}(x-y),(x-y)\rangle})}\,(1+\xi_y(x)),
\label{eq:G-asym1}
\end{equation}
where
$M_{1,\infty}(\xi_y,r;y)\leq c\,\max(\om(r),\sigma(r))$ for $0<r<\e$ 
 and $J_y(r)=I_y (r)-I_y(0)$ can also be written as
\begin{equation}
\frac{-1}{|S^{n-1}|}\int_{0<|z-y|<r}
\left({\rm tr}({\bf A}_z{\bf A}_y^{-1})-
n\,\frac{\langle {\bf A}_z {\bf A}_y^{-1/2}(z-y), {\bf A}_y^{-1/2}(z-y)\rangle}{|z-y|^2}\right)\,\frac{d z}{|z-y|^{n}}.
\label{eq:Jy(r)}
\end{equation}
However, if we choose coordinates in which $y=0$ and $a_{ij}(0)=\de_{ij}$, then we may apply Lemma 1 to absorb the exponential term in (\ref{eq:G-asym1}) into $1+\xi_y$ and obtain
\begin{equation}
G(x,y)=\frac{\langle{\bf A}_y^{-1}(x-y),(x-y)\rangle^{\frac{2-n}{2}}}{(n-2)|S^{n-1}|\,\sqrt{\det{\bf A}_y}}\,(1+\xi_y(x)) \quad\hbox{as}\ x\to y,
\label{eq:G-asym2}
\end{equation}
where
$M_{2,p}(\xi_y,r;y)\leq c_y\,\max(\om(r),\sigma(r),\theta(r))$ for $0<r<\e$
with $c_y$ depending on $\om$, $\theta$, $n$, and $I_y(0)$. But if we select a compact subset $K\subset U$, then $c_y$ may be taken independent of $y\in K$,
so we can replace $\xi_y(x)$ by $H(x,y)$ and obtain (\ref{eq:G-asym0}), (\ref{eq:M2p(H)}) as $|x-y|\to 0$.

Letting $F(x,y)=G(x,y)$ for $x,y\in U\subset V$, we have our fundamental solution in $U$.
\hfill $\Box$


\bigskip\noindent
{\bf Acknowledgements}
V.Maz'ya was supported by the UK
Engineering and Physical Sciences Research Council via the
research grant EP/F005563/1.

\end{document}